\documentclass[12pt]{amsart}
\usepackage{a4wide,hyperref,tikz}
\usepackage[english]{babel}

\newtheorem{theorem}{Theorem}
\newtheorem{proposition}{Proposition}
\newtheorem{lemma}{Lemma}
\newtheorem{corollary}{Corollary}
\theoremstyle{remark}
\newtheorem{example}{Example}

\title{On the structure of ($-\beta$)-integers}

\author{Wolfgang Steiner}
\address{LIAFA, CNRS, Universit\'e Paris Diderot -- Paris 7, Case 7014, 75205 Paris Cedex 13, France}
\email{steiner@liafa.jussieu.fr}

\begin{document}
\begin{abstract}
The $(-\beta)$-integers are natural generalisations of the $\beta$-integers, and thus of the integers, for negative real bases. 
When $\beta$ is the analogue of a Parry number, we describe the structure of the set of $(-\beta)$-integers by a fixed point of an anti-morphism.
\end{abstract}

\maketitle

\section{Introduction}
The aim of this paper is to study the structure of the set of real numbers having a digital expansion of the form 
\[
\sum_{k=0}^{n-1} a_k\, (-\beta)^k\,,
\]
where $(-\beta)$ is a negative real base with $\beta > 1$,  the digits $a_k \in \mathbb{Z}$ satisfy certain conditions specified below, and $n \ge 0$.
These numbers are called \emph{($-\beta$)-integers}, and have been recently studied by Ambro\v{z}, Dombek, Mas\'akov\'a and Pelantov\'a~\cite{Ambroz-Dombek-Masakova-Pelantova}.

Before dealing with these numbers, we recall some facts about \emph{$\beta$-integers}, which are the real numbers of the form 
\[
\pm \sum_{k=0}^{n-1} a_k\, \beta^k \quad \mbox{such that} \quad 0 \le \sum_{k=0}^{m-1} a_k\, \beta^k < \beta^m \quad \mbox{for all}\ 1 \le m \le n\,,
\]
i.e., $\sum_{k=0}^{n-1} a_k\, \beta^k$ is a greedy \emph{$\beta$-expansion}.
Equivalently, we can define the set of $\beta$-integers~as
\[
\mathbb{Z}_\beta = \mathbb{Z}_\beta^+ \cup (-\mathbb{Z}_\beta^+) \quad \mbox{with} \quad \mathbb{Z}_\beta^+ = \bigcup_{n\ge0} \beta^n\, T_\beta^{-n}(0)\,,
\]
where $T_\beta$ is the \emph{$\beta$-transformation}, defined by
\[
T_\beta:\,[0,1) \to [0,1)\,,\quad x \mapsto \beta x - \lfloor \beta x \rfloor\,.
\]
This map and the corresponding $\beta$-expansions were first studied by R\'enyi~\cite{Renyi57}.

The notion of $\beta$-integers was introduced in the domain of quasicrystallography, see for instance~\cite{Burdik-Frougny-Gazeau-Krejcar98}, and the structure of the $\beta$-integers is very well understood now.
We have $\beta\, \mathbb{Z}_\beta \subseteq \mathbb{Z}_\beta$, the set of distances between consecutive elements of $\mathbb{Z}_\beta$ is
\[
\Delta_\beta = \{T_\beta^n(1^-) \mid n \ge 0\}\,,
\]
where $T_\beta^n(x^-) = \lim_{y\to x-} T_\beta^n(y)$, and the sequence of distances between consecutive elements of $\mathbb{Z}_\beta^+$ is coded by the fixed point of a substition, see~\cite{Fabre95} for the case when $\Delta_\beta$ is a finite set, that is when $\beta$ is a \emph{Parry number}.
We give short proofs of these facts in Section~\ref{sec:beta-integers}.
More detailed properties of this sequence can be found e.g.\ in \cite{Balkova-Gazeau-Pelantova08,Balkova-Pelantova-Steiner08,Bernat-Masakova-Pelantova07,Frougny-Masakova-Pelantova04,Klouda-Pelantova09}.

Closely related to $\mathbb{Z}_\beta^+$ are the sets
\[
S_\beta(x) = \bigcup_{n\ge0} \beta^n\, T_\beta^{-n}(x) \qquad (x \in [0,1)),
\]
which were used by Thurston~\cite{Thurston89} to define (fractal) tilings of $\mathbb{R}^{d-1}$ when $\beta$ is a Pisot number of degree~$d$, i.e., a root $>1$ of a polynomial $x^d + p_1 x^{d-1} + \cdots + p_d \in \mathbb{Z}[x]$ such that all other roots have modulus~$<1$, and an algebraic unit, i.e., $p_d = \pm 1$.
These tilings allow e.g.\ to determine the $k$-th digit $a_k$ of a number without knowing the other digits, see~\cite{Kalle-Steiner}.

\begin{figure}[ht]
\centering
\begin{tikzpicture}[scale=4]
\draw(-.6667,0)node[left]{\small$0$}--(.3333,0) (0,-.6667)node[below]{\small$0\vphantom{/}$}--(0,.3333)
(-.6667,-.6667)node[below]{\small$-2/3$}node[left]{\small$-2/3$}--(.3333,-.6667)node[below]{\small$1/3$}--(.3333,.3333)--(-.6667,.3333)node[left]{\small$1/3$}--cycle
(-.1667,-.6667)node[below]{\small$-1/6$}--(-.1667,.3333);
\draw[very thick](-.6667,.3333)--(-.1667,-.6667) (-.1667,.3333)--(.3333,-.6667);

\begin{scope}[shift={(1.275,0)}]
\draw(-.618,0)node[left]{\small$0$}--(.382,0) (0,-.618)node[below]{\small$0\vphantom{/}$}--(0,.382)
(-.618,-.618)node[below]{\small$-1/\beta$}node[left]{\small$-1/\beta$}--(.382,-.618)node[below]{\small$1/\beta^2$}--(.382,.382)--(-.618,.382)node[left]{\small$1/\beta^2$}--cycle
(-.2361,-.618)node[below]{\small$-1/\beta^3$}--(-.2361,.382);
\draw[very thick](-.618,0)--(-.2361,-.618) (-.2361,.382)--(.382,-.618);
\end{scope}

\begin{scope}[shift={(2.55,0)}]
\draw(-.5698,0)node[left]{\small$0$}--(.4302,0) (0,-.5698)node[below]{\small$0\vphantom{\frac{1}{\beta+1}}$}--(0,.4302)
(-.5698,-.5698)node[below]{$\frac{-\beta}{\beta+1}$}node[left]{$\frac{-\beta}{\beta+1}$}--(.4302,-.5698)node[below]{$\frac{1}{\beta+1}$}--(.4302,.4302)--(-.5698,.4302)node[left=-.25em]{$\frac{1}{\beta+1}$}--cycle
(-.3247,-.5698)node[below=-.2ex]{$\frac{-\beta^{-1}}{\beta+1}$}--(-.3247,.4302);
\draw[very thick](-.5698,-.2451)--(-.3247,-.5698) (-.3247,.4302)--(.4302,-.5698);
\end{scope}
\end{tikzpicture}
\caption{The $(-\beta)$-transformation for $\beta = 2$ (left), $\beta = \frac{1+\sqrt5}{2} \approx 1.618$ (middle), and $\beta = \frac{1}{\beta} + \frac{1}{\beta^2} \approx 1.325$ (right).}
\label{f:Tmbeta}
\end{figure}
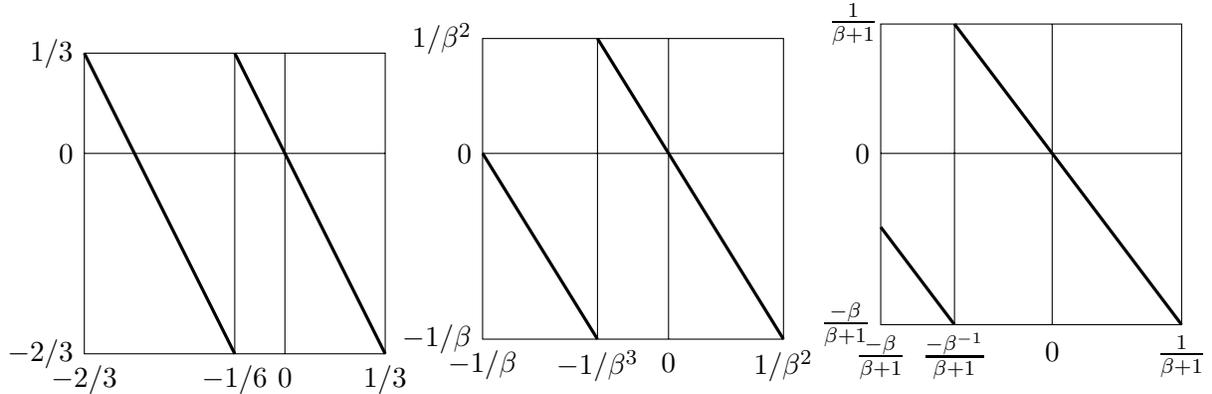

It is widely agreed that the greedy $\beta$-expansions are the natural representations of real numbers in a real base $\beta > 1$. 
For the case of negative bases, the situation is not so clear. 
Ito and Sadahiro~\cite{Ito-Sadahiro09} proposed recently to use the \emph{$(-\beta)$-transformation} defined by
\[
T_{-\beta}:\ \big[\tfrac{-\beta}{\beta+1}, \tfrac{1}{\beta+1}\big),\ x \mapsto -\beta x - \big\lfloor \tfrac{\beta}{\beta+1} -\beta x\big\rfloor\,,
\]
see also~\cite{Frougny-Lai09}.
This transformation has the important property that $T_{-\beta}(-x/\beta) = x$ for all $x \in \big(\frac{-\beta}{\beta+1}, \frac{1}{\beta+1}\big)$.
Some instances are depicted in Figures~\ref{f:Tmbeta}, \ref{f:Tgm2}, \ref{f:Tcomplex} and~\ref{f:Tcomplex2}.

The set of $(-\beta)$-integers is therefore defined by
\[
\mathbb{Z}_{-\beta} = \bigcup_{n\ge0} (-\beta)^n\, T_{-\beta}^{-n}(0)\,.
\]
These are the numbers
\[
\sum_{k=0}^{n-1} a_k\, (-\beta)^k \quad \mbox{such that} \quad  \frac{-\beta}{\beta+1} \le \sum_{k=0}^{m-1} a_k\, (-\beta)^{k-m} < \frac{1}{\beta+1} \quad \mbox{for all}\ 1 \le m \le n\,.
\]
Note that, in the case of $\beta$-integers, we have to add $-\mathbb{Z}_\beta^+$ to $\mathbb{Z}_\beta^+$ in order to obtain a set resembling~$\mathbb{Z}$.
In the case of $(-\beta)$-integers, this is not necessary because the $(-\beta)$-trans\-for\-mation allows to represent positive and negative numbers.

It is not difficult to see that $\mathbb{Z}_{-\beta} = \mathbb{Z} = \mathbb{Z}_\beta$ when $\beta \in \mathbb{Z}$, $\beta \ge 2$.
Some other properties of $\mathbb{Z}_{-\beta}$ can be found in~\cite{Ambroz-Dombek-Masakova-Pelantova}, mainly for the case when $T_{-\beta}^n\big(\frac{-\beta}{\beta+1}\big) \le 0$ and $T_{-\beta}^{2n-1}\big(\frac{-\beta}{\beta+1}\big) \ge \frac{1-\lfloor\beta\rfloor}{\beta}$ for all $n \ge 1$.
(Note that $T_{-\beta}^n\big(\frac{-\beta}{\beta+1}\big) \in \big(\frac{1}{\beta+1} - \frac{\lfloor\beta\rfloor}{\beta}, \frac{1-\lfloor\beta\rfloor}{\beta}\big) \cup \big(\frac{-\beta^{-1}}{\beta+1}, 0\big)$ implies $T_{-\beta}^{n+1}\big(\frac{-\beta}{\beta+1}\big) > 0$.)

The set 
\[
V_\beta = \big\{T_{-\beta}^n\big(\tfrac{-\beta}{\beta+1}\big) \mid n \ge 0\big\}
\]
plays a similar role for $(-\beta)$-expansions as the set $\{T_\beta^n(1^-) \mid n \ge 0\}$ for $\beta$-expansions.
If $V_\beta$ is a finite set, then we call $\beta > 1$ an \emph{Yrrap number}.
Note that these numbers are called \emph{Ito--Sadahiro numbers} in~\cite{Masakova-Pelantova}, in reference to~\cite{Ito-Sadahiro09}.
However, as the generalised $\beta$-trans\-for\-mations in \cite{Gora07} with $E = (1,\ldots,1)$ are, up to conjugation by the map $x \mapsto \frac{1}{\beta+1} - x$, the same as our $(-\beta)$-transformations, these numbers were already considered by G\'ora and perhaps by other authors.
Therefore, the neutral but intricate name \emph{$(-\beta)$-numbers} was chosen in~\cite{Liao-Steiner}, referring to the original name $\beta$-numbers for Parry numbers~\cite{Parry60}. 
The name Yrrap number, used in the present paper, refers to the connection with Parry numbers and to the fact that $T_{-\beta}$ is (locally) orientation-reversing. 

For any Yrrap number $\beta \ge (1+\sqrt5)/2$, we describe the sequence of $(-\beta)$-integers in terms of a two-sided infinite word  on a finite alphabet which is a fixed point of an anti-morphism (Theorem~\ref{t:main}).
Note that the orientation-reversing property of the map $x \mapsto -\beta x$ imposes the use of an anti-morphism instead of a morphism, and that anti-morphisms were considered in a similar context by Enomoto~\cite{Enomoto08}.

For $1 < \beta < \frac{1+\sqrt5}{2}$, we have $\mathbb{Z}_{-\beta} = \{0\}$, as already proved in~\cite{Ambroz-Dombek-Masakova-Pelantova}.
However, our study still makes sense for these bases $(-\beta)$ because we can also describe the sets 
\[
S_{-\beta}(x) = \lim_{n\to\infty} (-\beta)^n\, T_{-\beta}^{-n}(x) \qquad \big(x \in \big[\tfrac{-\beta}{\beta+1}, \tfrac{1}{\beta+1}\big)\big),
\]
where the limit set consists of the numbers lying in all but finitely many sets $(-\beta)^n\, T_{-\beta}^{-n}(x)$, $n \ge 0$.
Taking the limit instead of the union over all $n \ge 0$ implies that every $y \in \mathbb{R}$ lies in exactly one set $S_{-\beta}(x)$, $x \in \big[\tfrac{-\beta}{\beta+1}, \tfrac{1}{\beta+1}\big)$, see Lemma~\ref{l:Smbeta}.
Note that $T_{-\beta}^2\big(\frac{-\beta^{-1}}{\beta+1}\big) \ne \frac{-\beta}{\beta+1}$ when $\beta \not\in \mathbb{Z}$.
Other properties of the $(-\beta)$-transformation for $1 < \beta < \frac{1+\sqrt5}{2}$ are exhibited in~\cite{Liao-Steiner}.

\section{$\beta$-integers} \label{sec:beta-integers}
In this section, we consider the structure of $\beta$-integers. 
The results are not new, but it is useful to state and prove them in order to understand the differences with $(-\beta)$-integers.

Recall that $\Delta_\beta = \{T_\beta^n(1^-) \mid n \ge 0\}$, and let $\Delta_\beta^*$ be the free monoid generated by~$\Delta_\beta$.
Elements of $\Delta_\beta^*$ will be considered as words on the alphabet~$\Delta_\beta$, and the operation is the concatenation of words.
The \emph{$\beta$-substitution} is the morphism $\varphi_\beta:\, \Delta_\beta^* \to \Delta_\beta^*$, defined by 
\[
\varphi_\beta(x) = \!\!\! \underbrace{1 1 \cdots 1}_{\lceil\beta x\rceil-1\,\mathrm{times}} \!\!\! T_\beta(x^-) \qquad (x \in \Delta_\beta).
\]
Here, $1$~is an element of~$\Delta_\beta$ and not the identity element of~$\Delta_\beta^*$ (which is the empty word). 
Recall that, as $\varphi_\beta$ is a morphism, we have $\varphi_\beta(u v) = \varphi_\beta(u) \varphi_\beta(v)$ for all $u, v \in \Delta_\beta^*$.
Since $\varphi_\beta^{n+1}(1) = \varphi_\beta^n(\varphi_\beta(1))$ and $\varphi_\beta(1)$ starts with~$1$, $\varphi_\beta^n(1)$~is a prefix of~$\varphi_\beta^{n+1}(1)$ for every $n \ge 0$.

\begin{theorem} \label{t:beta}
For any $\beta > 1$, the set of non-negative $\beta$-integers takes the form
\[
\mathbb{Z}_\beta^+ = \{z_k \mid k \ge 0\} \quad \mbox{with} \quad z_k = \sum_{j=1}^k u_j\,,
\]
where $u_1 u_2 \cdots$ is the infinite word with letters in $\Delta_\beta$ which has $\varphi_\beta^n(1)$ as prefix for all $n \ge 0$.

The set of differences between consecutive $\beta$-integers is~$\Delta_\beta$.
\end{theorem}

The main observation for the proof of the theorem is the following.
We use the notation $|v| = k$ and $L(v) = \sum_{j=1}^k v_j$ for any $v = v_1 \cdots v_k \in \Delta_\beta^k$, $k \ge 0$.

\begin{lemma} \label{l:beta}
For any $n \ge 0$, $1 \le k \le |\varphi_\beta^n(1)|$, we have
\[
T_\beta^n\Big(\Big[\frac{z_{k-1}}{\beta^n}, \frac{z_k}{\beta^n}\Big)\Big) = [0, u_k)\,,
\]
and $z_{|\varphi_\beta^n(1)|} = L(\varphi_\beta^n(1)) = \beta^n$.  
\end{lemma}

\begin{proof}
For $n=0$, we have $|\varphi_\beta^0(1)| = 1$, $z_0 = 0$, $z_1 = 1$, $u_1 = 1$,  thus the statements are true.
Suppose that they hold for~$n$, and consider 
\[
u_1 u_2 \cdots u_{|\varphi_\beta^{n+1}(1)|} = \varphi_\beta^{n+1}(1) = \varphi_\beta(\varphi_\beta^n(1)) = \varphi_\beta(u_1)\, \varphi_\beta(u_2)\, \cdots\, \varphi_\beta(u_{|\varphi_\beta^n(1)|})\,.
\]
Let $1 \le k \le |\varphi_\beta^{n+1}(1)|$, and write $u_1 \cdots u_k = \varphi_\beta(u_1 \cdots u_{j-1})\, v_1 \cdots v_i$ with $1 \le j \le |\varphi_\beta^n(1)|$, $1 \le i \le |\varphi_\beta(u_j)|$, i.e., $v_1 \cdots v_i$ is a non-empty prefix of~$\varphi_\beta(u_j)$.

For any $x \in (0,1]$, we have $T_\beta(x^-) = \beta x - \lceil\beta x\rceil + 1$, hence $L(\varphi_\beta(x)) = \beta x$ for $x \in \Delta_\beta$.
This yields that 
\[
z_k = L(u_1 \cdots u_k) = \beta\, L(u_1 \cdots u_{j-1}) + L(v_1 \cdots v_i) = \beta z_{j-1} + i - 1 + v_i
\]
and $z_{k-1} = \beta z_{j-1} + i - 1$, hence 
\begin{gather*}
\Big[\frac{z_{k-1}}{\beta}, \frac{z_k}{\beta}\Big) = \Big[z_{j-1} + \frac{i-1}{\beta}, z_{j-1} + \frac{i-1+v_i}{\beta}\Big) \subseteq [z_{j-1}, z_{j-1} + u_j) = [z_{j-1}, z_j)\,, \\
T_\beta^{n+1}\Big(\Big[\frac{z_{k-1}}{\beta^{n+1}}, \frac{z_k}{\beta^{n+1}}\Big)\Big) = T_\beta\Big(\Big[\frac{i-1}{\beta}, \frac{i-1+v_i}{\beta}\Big)\Big) = [0, v_i) = [0, u_k)\,.
\end{gather*}
Moreover, we have $L(\varphi_\beta^{n+1}(1)) = \beta\, L(\varphi_\beta^n(1)) = \beta^{n+1}$, thus the statements hold for~$n+1$.
\end{proof}

\begin{proof}[Proof of Theorem~\ref{t:beta}]
By Lemma~\ref{l:beta}, we have $z_{|\varphi_\beta^n(1)|} = \beta^n$ for all $n \ge 0$, thus $[0,1)$ is split into the intervals $[z_{k-1}/\beta^n, z_k/\beta^n)$, $1 \le k \le |\varphi_\beta^n(1)|$.
Therefore, Lemma~\ref{l:beta}  yields that
\[
T_\beta^{-n}(0) = \{z_{k-1}/\beta^n \mid 1 \le k \le |\varphi_\beta^n(1)|\}\,,
\]
hence 
\[
\bigcup_{n\ge0} \beta^n\, T_\beta^{-n}(0) = \{z_k \mid k \ge 0\}\,.
\]
Since $u_k \in \Delta_\beta$ for all $k \ge 1$ and $u_{|\varphi^n(1)|} = T_\beta^n(1^-)$ for all $n \ge 0$, we have 
\[
\{z_k - z_{k-1} \mid k \ge 1\} = \{u_k \mid k \ge 1\} = \Delta_\beta\,. \qedhere
\]
\end{proof}

For the sets $S_\beta(x)$, Lemma~\ref{l:beta} gives the following corollary.

\begin{corollary} \label{c:Sbeta}
For any $x \in [0,1)$, we have 
\[
S_\beta(x) = \{z_k + x \mid k \ge 0,\, u_{k+1} > x\} \subseteq x + S_\beta(0)\,.
\]
\end{corollary}

In particular, we have $S_\beta(x) - x = S_\beta(y) - y$ for all $x, y \in [0,1)$ with $(x,y] \cap \Delta_\beta = \emptyset$.
From the definition of $S_\beta(x)$ and since $x \in \beta\, T_\beta^{-1}(x)$, we also get that
\[
S_\beta(x) = \bigcup_{y \in T_\beta^{-1}(x)} \beta\, S_\beta(y) \qquad \big(x \in [0,1)\big).
\]
This shows that $S_\beta(x)$ is the solution of a graph-directed iterated function system (GIFS) when $\beta$ is a Parry number, cf.~\cite[Section~3.2]{Kalle-Steiner}.

\section{$(-\beta)$-integers} \label{sec:-beta-integers}
We now turn to the discussion of $(-\beta)$-integers and sets $S_{-\beta}(x)$, $x \in \big[\frac{-\beta}{\beta+1}, \frac{1}{\beta+1}\big)$.

\begin{lemma} \label{l:Smbeta}
For any $\beta > 1$, $x \in \big[\frac{-\beta}{\beta+1}, \frac{1}{\beta+1}\big)$, we have
\[
S_{-\beta}(x) = \bigcup_{n \ge 0} (-\beta)^n\, \big(T_{-\beta}^{-n}(x) \setminus \big\{\tfrac{-\beta}{\beta+1}\big\}\big) = \bigcup_{y \in T_{-\beta}^{-1}(x)} (-\beta)\, S_{-\beta}(y)\,.
\]

For any $y \in \mathbb{R}$, there exists a unique $x \in \big[\frac{-\beta}{\beta+1}, \frac{1}{\beta+1}\big)$ such that $y \in S_{-\beta}(x)$. 

If $T_{-\beta}(x) = x$, then $S_{-\beta}(x) = \bigcup_{n \ge 0} (-\beta)^n\, T_{-\beta}^{-n}(x)$, in particular $S_{-\beta}(0) = \mathbb{Z}_{-\beta}$.
\end{lemma}

\begin{proof}
If $y \in S_{-\beta}(x)$, then we have $\frac{y}{(-\beta)^n} \in T_{-\beta}^{-n}(x)$ for all sufficiently large~$n$, thus $y \in (-\beta)^n\, \big(T_{-\beta}^{-n}(x) \setminus \big\{\tfrac{-\beta}{\beta+1}\big\}\big)$ for some $n \ge 0$.
On the other hand, $y \in (-\beta)^n\, \big(T_{-\beta}^{-n}(x) \setminus \big\{\tfrac{-\beta}{\beta+1}\big\}\big)$ for some $n \ge 0$ implies that $T_{-\beta}^m(\frac{y}{(-\beta)^m}) = T_{-\beta}^n(\frac{y}{(-\beta)^n}) = x$ for all $m \ge n$, thus $y \in S_{-\beta}(x)$. 
This shows the first equation.
Since $x \in \big(\frac{-\beta}{\beta+1}, \frac{1}{\beta+1}\big)$ implies that $x \in (-\beta)\, \big(T_{-\beta}^{-1}(x) \setminus \big\{\tfrac{-\beta}{\beta+1}\big\}\big)$, we obtain that $S_{-\beta}(x) = \bigcup_{y \in T_{-\beta}^{-1}(x)} (-\beta)\, S_{-\beta}(y)$ for all $x \in \big[\frac{-\beta}{\beta+1}, \frac{1}{\beta+1}\big)$.

For any $y \in \mathbb{R}$, we have $y \in S_{-\beta}\big(T_{-\beta}^n\big(\frac{y}{(-\beta)^n}\big)\big)$ for all $n \ge 0$ such that $\frac{y}{(-\beta)^n} \in \big(\frac{-\beta}{\beta+1}, \frac{1}{\beta+1}\big)$, thus $y \in S_{-\beta}(x)$ for some $x \in \big[\frac{-\beta}{\beta+1}, \frac{1}{\beta+1}\big)$.
To show that this $x$ is unique, let $y \in S_{-\beta}(x)$ and $y \in S_{-\beta}(x')$ for some $x, x' \in \big[\frac{-\beta}{\beta+1}, \frac{1}{\beta+1}\big)$.
Then we have $y \in (-\beta)^n\, \big(T_{-\beta}^{-n}(x) \setminus \big\{\tfrac{-\beta}{\beta+1}\big\}\big)$ and $y \in (-\beta)^m\, \big(T_{-\beta}^{-m}(x') \setminus \big\{\tfrac{-\beta}{\beta+1}\big\}\big)$ for some $m, n \ge 0$, thus $x = T_{-\beta}^n\big(\frac{y}{(-\beta)^n}\big) = T_{-\beta}^m\big(\frac{y}{(-\beta)^m}\big) = x'$.

If $T_{-\beta}^n\big(\frac{-\beta}{\beta+1}\big) = x = T_{-\beta}(x)$, then $T_{-\beta}^{n+2}\big(\tfrac{-\beta^{-1}}{\beta+1}\big) = T_{-\beta}^{n+1}\big(\tfrac{-\beta}{\beta+1}\big) = T_{-\beta}(x) = x$ yields that $(-\beta)^n\, \frac{-\beta}{\beta+1} \in S_{-\beta}(x)$, which shows that $S_{-\beta}(x) = \bigcup_{n \ge 0} (-\beta)^n\, T_{-\beta}^{-n}(x)$ when $T_{-\beta}(x) = x$.
\end{proof}

The first two statements of the following proposition can also be found in~\cite{Ambroz-Dombek-Masakova-Pelantova}.

\begin{proposition} \label{p:mbeta}
For any $\beta > 1$, we have $(-\beta)\, \mathbb{Z}_{-\beta} \subseteq \mathbb{Z}_{-\beta}$.

If $\beta < (1+\sqrt5)/2$, then $\mathbb{Z}_{-\beta} = \{0\}$.

If $\beta \ge (1+\sqrt5)/2$, then 
\[
\mathbb{Z}_{-\beta} \cap (-\beta)^n\, [-\beta, 1] = \big\{(-\beta)^n, (-\beta)^{n+1}\big\} \cup (-\beta)^{n+2}\, \big(T_{-\beta}^{-n-2}(0) \cap \big(\tfrac{-1}{\beta}, \tfrac{1}{\beta^2}\big)\big)
\]
for all $n \ge 0$, in particular
\[
\mathbb{Z}_{-\beta} \cap [-\beta, 1]  = \left\{\begin{array}{cl}\{-\beta, -\beta + 1, \ldots, -\beta + \lfloor\beta\rfloor, 0, 1\} & \mbox{if}\ \beta^2 \ge \lfloor\beta\rfloor (\beta + 1), \\[1ex] \{-\beta, -\beta + 1, \ldots, -\beta + \lfloor\beta\rfloor - 1, 0, 1\} & \mbox{if}\ \beta^2 < \lfloor\beta\rfloor (\beta + 1).\end{array}\right.
\]
\end{proposition}

\begin{proof}
The inclusion $(-\beta)\, \mathbb{Z}_{-\beta} \subseteq \mathbb{Z}_{-\beta}$ is a consequence of Lemma~\ref{l:Smbeta} and $0 \in T_{-\beta}^{-1}(0)$. 

If $\beta < \frac{1+\sqrt5}{2}$, then $\frac{-1}{\beta} < \frac{-\beta}{\beta+1}$, hence $T_{-\beta}^{-1}(0) = \{0\}$ and $\mathbb{Z}_{-\beta} = \{0\}$, see Figure~\ref{f:Tmbeta} (right).

If $\beta \ge \frac{1+\sqrt5}{2}$, then $\frac{-1}{\beta} \in T_{-\beta}^{-1}(0)$ implies $1 \in \mathbb{Z}_{-\beta}$, thus $(-\beta)^n \in \mathbb{Z}_{-\beta}$ for all $n \ge 0$.
Clearly,
\[
(-\beta)^{n+2}\, \big(T_{-\beta}^{-n-2}(0) \cap \big(\tfrac{-1}{\beta}, \tfrac{1}{\beta^2}\big)\big) \subseteq \mathbb{Z}_{-\beta} \cap (-\beta)^n\, (-\beta, 1)\,.
\]
To show the other inclusion, let $z \in (-\beta)^m\, T_{-\beta}^{-m}(0) \cap (-\beta)^n\, (-\beta, 1)$ for some $m \ge 0$.
If $z \ne (-\beta)^m \frac{-\beta}{\beta+1}$, then $\frac{z}{(-\beta)^m} \in \big(\frac{-\beta}{\beta+1}, \frac{1}{\beta+1}\big)$ and $\frac{z}{(-\beta)^{n+2}} \in \big(\frac{-1}{\beta}, \frac{1}{\beta^2}\big) \subseteq \big(\frac{-\beta}{\beta+1}, \frac{1}{\beta+1}\big)$ imply that $T_{-\beta}^{n+2}\big(\frac{z}{(-\beta)^{n+2}}\big) = 
T_{-\beta}^m\big(\frac{z}{(-\beta)^m}\big) = 0$.
If $z = (-\beta)^m \frac{-\beta}{\beta+1}$, then 
\[
T_{-\beta}^{n+2}\big(\tfrac{z}{(-\beta)^{n+2}}\big) =  T_{-\beta}^{n+2}\big(\tfrac{(-\beta)^{m-n-1}}{\beta+1}\big) = T_{-\beta}^{m+2}\big(\tfrac{-\beta^{-1}}{\beta+1}\big) = T_{-\beta}^{m+1}\big(\tfrac{-\beta}{\beta+1}\big) = T_{-\beta}(0) = 0\,,
\]
where we have used that $\frac{z}{(-\beta)^{n+2}} \in \big(\frac{-\beta}{\beta+1}, \frac{1}{\beta+1}\big)$ implies $m \le n$.
Therefore, we have $z \in (-\beta)^{n+2}\, T_{-\beta}^{-n-2}(0)$ for all $z \in \mathbb{Z}_{-\beta} \cap (-\beta)^n\, (-\beta, 1)$.

Consider now $n = 0$, then 
\[
\mathbb{Z}_{-\beta} \cap [-\beta, 1] = \{-\beta, 1\} \cup \{z \in (-\beta, 1) \mid T_{-\beta}^2(z/\beta^2) = 0\}\,.
\]
Since $\frac{-\lfloor\beta\rfloor}{\beta} \ge \frac{-\beta}{\beta+1}$ if and only if $\beta^2 \ge \lfloor\beta\rfloor (\beta + 1)$, we obtain that
\[
(-\beta)\, T_{-\beta}^{-1}(0)  = \left\{\begin{array}{cl}\{0, 1, \ldots, \lfloor\beta\rfloor\} & \mbox{if}\ \beta^2 \ge \lfloor\beta\rfloor (\beta + 1), \\[1ex] \{0, 1, \ldots, \lfloor\beta\rfloor - 1\} & \mbox{if}\ \beta^2 < \lfloor\beta\rfloor (\beta + 1).\end{array}\right.
\]
If $T_{-\beta}^2(z/\beta^2) = 0$, then $z = - a_1 \beta + a_0$ with $a_0 \in (-\beta)\, T_{-\beta}^{-1}(0)$, $a_1 \in \{0, 1, \ldots, \lfloor\beta\rfloor\}$, and $\mathbb{Z}_{-\beta} \cap [-\beta, 1]$ consists of those numbers $- a_1 \beta + a_0$ lying in $[-\beta, 1]$.
\end{proof}

Proposition~\ref{p:mbeta} shows that the maximal difference between consecutive $(-\beta)$-integers exceeds~$1$ whenever $\beta^2 < \lfloor\beta\rfloor (\beta + 1)$, i.e., $T_{-\beta}\big(\frac{-\beta}{\beta+1}\big) < 0$.
For an example, this was also proved in~\cite{Ambroz-Dombek-Masakova-Pelantova}.
In Examples~\ref{x:complex} and~\ref{x:complex2}, we see that the distance between two consecutive $(-\beta)$-integers can be even greater than~$2$, and the structure of $\mathbb{Z}_{-\beta}$ can be quite complicated.
Therefore, we adapt a slightly different strategy as for~$\mathbb{Z}_\beta$.

\medskip
In the following, we always assume that the set
\[
V'_\beta = V_\beta \cup \{0\} = \big\{T_{-\beta}^n\big(\tfrac{-\beta}{\beta+1}\big) \mid n \ge 0\big\} \cup \{0\}
\]
is finite, i.e., $\beta$ is an Yrrap number, and let $\beta$ be fixed. 
For $x \in V'_\beta$, let 
\[
r_x = \min\big\{y \in V'_\beta \cup \big\{\tfrac{1}{\beta+1}\big\} \mid y > x\big\}\,, \quad \widehat{x} = \tfrac{x+r_x}{2}, \quad  J_x = \{x\} \quad \mbox{and} \quad  J_{\widehat{x}} = (x, r_x)\,.
\]
Then $\{J_a \mid a \in A_\beta\}$ forms a partition of $\big[\frac{-\beta}{\beta+1}, \frac{1}{\beta+1}\big)$, where
\[
A_\beta = V'_\beta \cup \,\widehat{\!V'_\beta\!}\,\,, \quad \mbox{with} \quad \,\widehat{\!V'_\beta\!}\, = \{\widehat{x} \mid x \in V'_\beta\}\,.
\]
We have $T_{-\beta}(J_x) = J_{T_{-\beta}(x)}$ for every $x \in V'_\beta$, and the following lemma shows that the image of every~$J_{\widehat{x}}$, $x \in V'_\beta$, is a union of intervals~$J_a$, $a \in A_\beta$.

\begin{lemma} \label{l:hatx}
Let $x \in V'_\beta$ and write
\[
J_{\widehat{x}} \cap T_{-\beta}^{-1}(V'_\beta) = \{y_1, \ldots, y_m\}\,, \quad \mbox{with} \quad x = y_0 < y_1 < \cdots < y_m < y_{m+1} = r_x\,.
\]
For any $0 \le i \le m$, we have
\[
T_{-\beta}\big((y_i, y_{i+1})\big) = J_{\,\widehat{\!x_i\!}\,} \quad \mbox{with} \quad x_i = \lim_{y\to y_{i+1}-} T_{-\beta}(y)\,, \ \mbox{i.e.},\ \,\widehat{\!x_i\!}\, = T_{-\beta}\big(\tfrac{y_i+y_{i+1}}{2}\big)\,,
\]
and $\beta (y_{i+1} - y_i) = \lambda(J_{\,\widehat{\!x_i\!}\,})$, where $\lambda$ denotes the Lebesgue measure.
\end{lemma}

\begin{proof}
Since $T_{-\beta}$ maps no point in $(y_i, y_{i+1})$ to $\frac{-\beta}{\beta+1} \in V'_\beta$, the map is continuous on this interval and $\lambda(T_{-\beta}((y_i, y_{i+1}))) = \beta (y_{i+1} - y_i)$.
We have $x_i \in V'_\beta$ since $x_i = T_{-\beta}(y_{i+1})$ in case $y_{i+1} < \frac{1}{\beta+1}$, and $x_i = \frac{-\beta}{\beta+1}$ in case $y_{i+1} = \frac{1}{\beta+1}$.
Since $y_i = \max\{y \in T_{-\beta}^{-1}(V'_\beta) \mid y < y_{i+1}\}$, we obtain that $r_{x_i}  = \lim_{y\to y_i+}T_{-\beta}(y)$, thus $T_{-\beta}((y_i, y_{i+1})) = (x_i, r_{x_i})$.
\end{proof}

In view of Lemma~\ref{l:hatx}, we define an anti-morphism $\psi_\beta:\, A_\beta^* \to A_\beta^*$ by
\[
\psi_\beta(x) = T_{-\beta}(x) \quad \mbox{and} \quad \psi_\beta(\widehat{x}) = \,\widehat{\!x_m\!}\, \, T_{-\beta}(y_m) \cdots \widehat{x_1} \, T_{-\beta}(y_1) \, \widehat{x_0} \qquad (x \in V'_\beta),
\]
with $m$, $x_i$ and $y_i$ as in Lemma~\ref{l:hatx}.
Here, anti-morphism means that $\psi_\beta(u v) = \psi_\beta(v) \psi_\beta(u)$ for all $u, v \in A_\beta^*$.
Now, the last letter of $\psi_\beta(\widehat{0})$ is~$\widehat{t}$, with $t = \max\{x \in V_\beta \mid x < 0\}$, and the first letter of $\psi_\beta(\widehat{t}\,)$ is~$\widehat{0}$.
Therefore, $\psi_\beta^{2n}(\widehat{0})$ is a prefix of $\psi_\beta^{2n+2}(\widehat{0}) = \psi_\beta^{2n}(\psi_\beta^2(\widehat{0}))$ and $\psi_\beta^{2n+1}(\widehat{0})$ is a suffix of $\psi_\beta^{2n+3}(\widehat{0})$ for every $n \ge 0$.

\begin{theorem} \label{t:mbeta}
For any Yrrap number $\beta \ge (1+\sqrt5)/2$, we have
\[
\mathbb{Z}_{-\beta} = \{z_k \mid k \in \mathbb{Z},\, u_{2k} = 0\} \quad \mbox{with} \quad z_k = \left\{\begin{array}{cl}\sum_{j=1}^k \lambda(J_{u_{2j-1}}) & \mbox{if}\ k \ge 0\,, \\[1ex] - \sum_{j=1}^{|k|} \lambda(J_{u_{-2j+1}}) & \mbox{if}\ k \le 0\,,\end{array}\right.
\]
where $\cdots u_{-1} u_0 u_1 \cdots$ is the two-sided infinite word on the finite alphabet $A_\beta$ such that $u_0 = 0$,  $\psi_\beta^{2n}(\widehat{0})$ is a prefix of $u_1 u_2 \cdots$ and $\psi_\beta^{2n+1}(\widehat{0})$ is a suffix of $\cdots u_{-2} u_{-1}$ for all $n \ge 0$. 
\end{theorem}

Note that $\cdots u_{-1} u_0 u_1 \cdots$ is a fixed point of~$\psi_\beta$, with $u_0$ being mapped to~$u_0$.

The following lemma is the analogue of Lemma~\ref{l:beta}.
We use the notation 
\[
L(v) = \sum_{j=1}^k \lambda(J_{v_j}) \quad \mbox{if} \ v = v_1 \cdots v_k \in A_\beta^k\,.
\]
Note that $u_{2k} \in V'_\beta$ and $u_{2k+1} \in \,\widehat{\!V'_\beta\!}\,$ for all $k \in \mathbb{Z}$, thus $\lambda(J_{u_{2k}}) = 0$ for all $k \in \mathbb{Z}$.

\begin{lemma} \label{l:mbeta}
For any $n \ge 0$, $0 \le k < |\psi_\beta^n(\widehat{0})|/2$, we have
\[
T_{-\beta}^n\Big(\frac{z_{(-1)^n k}}{(-\beta)^n}\Big) = u_{(-1)^n2k}\,, \quad
T_{-\beta}^n\Big(\Big(\frac{z_{(-1)^n k}}{(-\beta)^n}, \frac{z_{(-1)^n(k+1)}}{(-\beta)^n}\Big)\Big) = J_{u_{(-1)^n(2k+1)}}\,,
\]
and $z_{(-1)^n (|\psi_\beta^n(\widehat{0})|+1)/2} = (-1)^n\, L\big(\psi_\beta^n(\widehat{0})\big) = \lambda(J_{\widehat{0}})\, (-\beta)^n = r_0\, (-\beta)^n$. 
\end{lemma}

\begin{proof}
The statements are true for $n=0$ since $|\psi_\beta^0(\widehat{0})| = 1$, $z_0 = 0$, $z_1 = \lambda(J_{\widehat{0}}) = r_0$.

Suppose that they hold for even~$n$, and consider 
\[
u_{-|\psi_\beta^{n+1}(\widehat{0})|} \cdots u_{-2} u_{-1} = \psi_\beta^{n+1}(\widehat{0}) = \psi_\beta\big(\psi_\beta^n(\widehat{0})\big) = \psi_\beta(u_{|\psi_\beta^n(\widehat{0})|}) \cdots \psi_\beta(u_2) \psi_\beta(u_1)\,.
\]
Let $0 \le k < |\psi_\beta^{n+1}(\widehat{0})|/2$, and write 
\[
u_{-2k-1} \cdots u_{-1} = v_{-2i-1} \cdots v_{-1}\, \psi_\beta(u_1 \cdots u_{2j})
\]
with $0 \le j < |\psi_\beta^n(\widehat{0})|/2$, $0 \le i < |\psi_\beta(u_{2j+1})|/2$, i.e., $u_{-2i-1} \cdots u_{-1}$ is a suffix of~$\psi_\beta(u_{2j+1})$.

By Lemma~\ref{l:hatx}, we have $L(\psi_\beta(\widehat{x})) = \beta \, \lambda(J_{\widehat{x}})$ for any $x \in V'_\beta$.
This yields that 
\[
- z_{-k-1} = \beta\, L(u_1 \cdots u_{2j}) + L(v_{-2i-1} \cdots v_{-1}) = \beta\, z_j + L(v_{-2i-1} \cdots v_{-1})
\]
and $-z_{-k} = \beta\, z_j + L(v_{-2i} \cdots v_{-1})$.
By the induction hypothesis, we obtain that
\begin{align*}
T_{-\beta}^{n+1}\Big(\frac{z_{-k}}{(-\beta)^{n+1}}\Big) & = T_{-\beta}^{n+1}\bigg(\frac{z_j}{(-\beta)^n} - \frac{L(v_{-2i}\cdots v_{-1})}{(-\beta)^{n+1}}\bigg) \\
& = \left\{\begin{array}{ll}T_{-\beta}(u_{2j}) = \psi_{\beta}(u_{2j}) = u_{-2k} & \mbox{if}\ i = 0, \\[1ex] T_{-\beta}\big(x + L(v_{-2i} \cdots v_{-1})/\beta\big) = T_{-\beta}(y_i) = v_{-2i} = u_{-2k} & \mbox{if}\ i \ge 1,\end{array}\right.
\end{align*}
where the $y_i$'s are the numbers from Lemma~\ref{l:hatx} for $\widehat{x} = u_{2j+1}$, and
\[
T_{-\beta}^{n+1}\Big(\Big(\frac{z_{-k}}{(-\beta)^{n+1}}, \frac{z_{-k-1}}{(-\beta)^{n+1}}\Big)\Big) = T_{-\beta}\big((y_i, y_{i+1})\big) = J_{v_{-2i-1}}  = J_{u_{-2k-1}}\,.
\]
Moreover, we have $L(\psi_\beta^{n+1}(\widehat{0})) = \beta\, L(\psi_\beta^n(\widehat{0})) = r_0 \beta^{n+1}$, thus the statements hold for~$n+1$.

The proof for odd $n$ runs along the same lines and is therefore omitted.
\end{proof}

\begin{proof}[Proof of Theorem~\ref{t:mbeta}]
By Lemma~\ref{l:mbeta}, we have $z_{(-1)^n(|\psi_\beta^n(\widehat{0})|+1)/2} = r_0\, (-\beta)^n$ for all $n \ge 0$, thus $[0,r_0)$ splits into the intervals $\big(z_{(-1)^n k} (-\beta)^{-n}, z_{(-1)^n(k+1)} (-\beta)^{-n}\big)$ and points $z_{(-1)^n k} (-\beta)^{-n}$, $0 \le k < |\psi_\beta^n(\widehat{0})|/2$, hence
\[
T_{-\beta}^{-n}(0) \cap [0,r_0) = \big\{z_{(-1)^n k} (-\beta)^{-n} \mid 0 \le k < |\psi_\beta^n(\widehat{0})|/2,\, u_{(-1)^n2k} = 0\big\}\,.
\]
Let $m \ge 1$ be such that $\beta^{2m} r_0 \ge \frac{1}{\beta+1}$.
Then we have $\big(\frac{-\beta}{\beta+1}, \frac{1}{\beta+1}\big) \subseteq (-\beta^{2m+1} r_0, \beta^{2m} r_0)$, and 
\[
T_{-\beta}^{-n}(0) \setminus \big\{\tfrac{-\beta}{\beta+1}\big\} \subseteq (-\beta)^{2m} \big(T_{-\beta}^{-n-2m}(0) \cap [0,r_0)\big) \cup (-\beta)^{2m+1} \big(T_{-\beta}^{-n-2m-1}(0) \cap [0,r_0)\big)\,,
\]
thus
\[
\bigcup_{n\ge0} (-\beta)^n\, \big(T_{-\beta}^{-n}(0) \setminus \big\{\tfrac{-\beta}{\beta+1}\big\}\big) = \bigcup_{n\ge0} (-\beta)^n\, \big(T_{-\beta}^{-n}(0) \cap [0,r_0)\big) = \{z_k \mid k \in \mathbb{Z},\, u_{2k} = 0\}\,. 
\]
Together with Lemma~\ref{l:Smbeta}, this proves the theorem.
\end{proof}

As in the case of positive bases, the word $\cdots u_{-1} u_0 u_1 \cdots$ also describes the sets $S_{-\beta}(x)$.
Theorem~\ref{t:mbeta} and Lemma~\ref{l:mbeta} give the following corollary.

\begin{corollary} \label{c:Smbeta}
For any $x \in V'_\beta$, $y \in J_{\widehat{x}}$, we have 
\[
S_{-\beta}(x) = \{z_k \mid k \in \mathbb{Z},\, u_{2k} = x\} \quad \mbox{and} \quad S_{-\beta}(y) = \{z_k + y - x \mid k \in \mathbb{Z},\, u_{2k+1} = \widehat{x}\} \,.
\]
\end{corollary}

Lemma~\ref{l:Smbeta} and Corollary~\ref{c:Smbeta} imply that $S_{-\beta}(x)$ is the solution of a GIFS for any Yrrap number $\beta \ge (1+\sqrt5)/2$, $x \in \big[\frac{-\beta}{\beta+1}, \frac{1}{\beta+1}\big)$, cf.\ the end of Section~\ref{sec:beta-integers}.

Recall that our main goal is to understand the structure of $\mathbb{Z}_{-\beta}$ (in case $\beta \ge (1+\sqrt5)/2$), i.e., to describe the occurrences of $0$ in the word $\cdots u_{-1} u_0 u_1 \cdots$ defined in Theorem~\ref{t:mbeta} and the words between two successive occurrences.
Let
\[
R_\beta = \{u_{2k} u_{2k+1} \cdots u_{2s(k)-1} \mid k \in \mathbb{Z},\, u_{2k} = 0\} \quad \mbox{with} \quad s(k) = \min\{j \in \mathbb{Z} \mid u_{2j} = 0,\, j > k\}
\]
be the set of return words of~$0$ in $\cdots u_{-1} u_0 u_1 \cdots$.
Note that $s(k)$ is defined for all $k \in \mathbb{Z}$ since $(-\beta)^n \in \mathbb{Z}_{-\beta}$ for all $n \ge 0$ by Proposition~\ref{p:mbeta}.

For any $w \in R_\beta$, the word $\psi_\beta(w0)$ is a factor of $\cdots u_{-1} u_0 u_1 \cdots$ starting and ending with~$0$, thus we can write $\psi_\beta(w0) = w_1 \cdots w_m 0$ with $w_j \in R_\beta$, $1 \le j \le m$, and set
\[
\varphi_{-\beta}(w) = w_1 \cdots w_m\,.
\]
This defines an anti-morphism $\varphi_{-\beta}:\, R_\beta^* \to R_\beta^*$, which plays the role of the $\beta$-substitution.

Since $\cdots u_{-1} u_0 u_1 \cdots$ is generated by $u_1 = \widehat{0}$, as described in Theorem~\ref{t:mbeta}, we consider $w_\beta = u_0 u_1 \cdots u_{2s(0)-1}$.
We have
\[
[0, 1] = \big[0, \tfrac{1}{\beta+1}\big) \cup \big[\tfrac{1}{\beta+1}, 1\big]\,, \quad T_{-\beta}\big((-\beta)^{-1} \big[\tfrac{1}{\beta+1}, 1\big]\big) = \big[\tfrac{-\beta}{\beta+1}, 0\big]\,,
\]
thus $L(w_\beta) = 1$ and
\[
w_\beta = 0 \, \widehat{0} \, x_1 \, \widehat{x_1} \, \cdots \, x_m \, \,\widehat{\!x_m\!}\, \, x_{-\ell} \, \,\widehat{\!x_{-\ell}\!\!}\,\,  \, \cdots x_{-1} \, \,\,\widehat{\!\!x_{-1}\!\!}\,\,\,,
\]
where the $x_j$ are defined by $V'_\beta = \{x_{-\ell}, \ldots, x_{-1}, 0, x_1, \ldots, x_m\}$, $x_{-\ell} < \cdots < x_{-1} < 0 < x_1 < \cdots < x_m$.

\begin{theorem} \label{t:main}
For any Yrrap number $\beta \ge (1+\sqrt5)/2$, we have
\[
\mathbb{Z}_{-\beta} = \{z'_k \mid k \in \mathbb{Z}\} \quad \mbox{with} \quad z'_k = \left\{\begin{array}{cl}\sum_{j=1}^k L(u'_j) & \mbox{if}\ k \ge 0\,, \\[1ex] - \sum_{j=1}^{|k|} L(u'_{-j}) & \mbox{if}\ k \le 0\,,\end{array}\right. 
\]
where $\cdots u'_{-2} u'_{-1} u'_1 u'_2 \cdots$ is the two-sided infinite word on the finite alphabet $R_\beta$ such that $\varphi_{-\beta}^{2n}(w_\beta)$ is a prefix of $u'_1 u'_2 \cdots$ and $\varphi_{-\beta}^{2n+1}(w_\beta)$ is a suffix of $\cdots u'_{-2} u'_{-1}$ for all $n \ge 0$.

The set of distances between consecutive $(-\beta)$-integers is
\[
\Delta_{-\beta} = \{z'_{k+1} - z'_k \mid k \in \mathbb{Z}\} = \{L(w) \mid w \in R_\beta\}\,.
\]
\end{theorem}

Note that the index $0$ is omitted in $\cdots u'_{-2} u'_{-1} u'_1 u'_2 \cdots$ for reasons of symmetry.

\begin{proof}
The definitions of $\cdots u_{-1} u_0 u_1 \cdots$ in Theorem~\ref{t:mbeta} and of $\varphi_{-\beta}$ imply that $\cdots u'_{-2} u'_{-1}$ $u'_1 u'_2 \cdots$ is the derived word of $\cdots u_{-1} u_0 u_1 \cdots$ with respect to~$R_\beta$, that is
\[
u'_k = u_{|u'_1 \cdots u'_{k-1}|} \cdots u_{|u'_1 \cdots u'_k|-1}\,, \quad u'_{-k} = u_{-|u'_{-k} \cdots u'_{-1}|} \cdots u_{-|u'_{1-k} \cdots u'_{-1}|-1} \quad (k \ge 1)
\]
with 
\[
\{|u'_1 \cdots u'_{k-1}| \mid k \ge 1\} \cup \{-|u'_{-k} \cdots u'_{-1}| \mid k \ge 1\} = \{k \in \mathbb{Z} \mid u_k = 0\}\,.
\]
Here, for any $v \in R_\beta^*$, $|v|$~denotes the length of $v$ as a word in $A_\beta^*$, not in~$R_\beta^*$.
Since
\[
z'_k = \sum_{j=1}^k L(u'_j) = \hspace{-4mm} \sum_{j=0}^{|u'_1 \cdots u'_k|-1} \hspace{-4mm} \lambda(J_{u_j}) = \hspace{-2mm} \sum_{j=1}^{|u'_1 \cdots u'_k|} \hspace{-2mm} \lambda(J_{u_j})\,, \quad z'_{-k} = - \sum_{j=1}^k L(u'_{-j}) = - \hspace{-4mm} \sum_{j=1}^{|u'_{-k} \cdots u'_{-1}|} \hspace{-4mm} \lambda(J_{u_{-j}})
\]
for all $k \ge 0$, Theorem~\ref{t:mbeta} yields that $\{z'_k \mid k \in \mathbb{Z}\} = \mathbb{Z}_{-\beta}$.

It follows from the definition of $R_\beta$ that $\Delta_{-\beta} = \{L(w) \mid w \in R_\beta\}$.

It remains to show that $R_\beta$ is a finite set.
We first show that the restriction of $\psi_\beta$ to $\,\widehat{\!V'_\beta\!}\,$ is primitive, which means that there exists some $m \ge 1$ such that, for every $x \in V'_\beta$, $\psi_\beta^m(\widehat{x})$ contains all elements of $\,\widehat{\!V'_\beta\!}\,$.
The proof is taken from \cite[Proposition~8]{Gora07}.
If $\beta > 2$, then the largest connected pieces of images of $J_{\widehat{x}}$ under $T_{-\beta}$ grow until they cover two consecutive discontinuity points $\frac{1}{\beta+1} - \frac{a+1}{\beta}$, $\frac{1}{\beta+1} - \frac{a}{\beta}$ of $T_{-\beta}$, and the next image covers all intervals $J_{\widehat{y}}$, $y \in V'_\beta$.
If $\frac{1+\sqrt5}{2} < \beta \le 2$, then $\beta^2 > 2$ implies that the largest connected pieces of images of $J_{\widehat{x}}$ under $T_{-\beta}^2$ grow until they cover two consecutive discontinuity points of $T_{-\beta}^2$.
Since
\begin{align*}
T_{-\beta}^2\big(\big(\tfrac{-\beta}{\beta+1}, \tfrac{\beta^{-2}}{\beta+1} - \tfrac{1}{\beta}\big)\big) & = \big(\tfrac{-\beta^3+\beta^2+\beta}{\beta+1}, \tfrac{1}{\beta+1}\big)\,, &
T_{-\beta}^2\big(\big(\tfrac{\beta^{-2}}{\beta+1} - \tfrac{1}{\beta}, \tfrac{-\beta^{-1}}{\beta+1}\big)\big) & = \big(\tfrac{-\beta}{\beta+1},  \tfrac{\beta^2-\beta-1}{\beta+1}\big)\,, \\
T_{-\beta}^2\big(\big(\tfrac{-\beta^{-1}}{\beta+1}, \tfrac{\beta^{-2}}{\beta+1}\big)\big) & = \big(\tfrac{-\beta}{\beta+1}, \tfrac{1}{\beta+1}\big)\,,  & 
T_{-\beta}^2\big(\big(\tfrac{\beta^{-2}}{\beta+1}, \tfrac{1}{\beta+1}\big)\big) & = \big(\tfrac{-\beta}{\beta+1},  \tfrac{\beta^2-\beta-1}{\beta+1}\big)\,,  
\end{align*}
the next image covers the fixed point~$0$, and further images grow until after a finite number of steps they cover all intervals $J_{\widehat{y}}$, $y \in V'_\beta$.
The case $\beta = \frac{1+\sqrt5}{2}$ is treated in Example~\ref{x:gm}.

If $T_{-\beta}^n\big(\frac{-\beta}{\beta+1}\big) \ne 0$ for all $n \ge 0$, then $T_{-\beta}^n$ is continuous at all points $x \in \big(\tfrac{-\beta}{\beta+1},  \tfrac{1}{\beta+1}\big)$ with $T_{-\beta}^n(x) = 0$, thus $u_{2k} = 0$ is equivalent to $u_{2k+1} = \widehat{0}$ (see also Proposition~\ref{p:hat} below).
Hence we can consider the return words of $\widehat{0}$ in $\cdots u_{-1} u_0 u_1 \cdots$ instead of the return words of~$0$.
Since $\psi_\beta^m(\widehat{x_0}\, x_1 \, \widehat{x_2})$ has at least two occurrences of $\widehat{0}$ for all $x_0, x_1, x_2 \in V'_\beta$, there are only finitely many such return words.
If $T_{-\beta}^n\big(\frac{-\beta}{\beta+1}\big) = 0$, then $\psi_\beta^n(x_0\, \widehat{x_1} \, x_2)$ starts and ends with~$0$ for all $x_0, x_1, x_2 \in V'_\beta$, hence $R_\beta$ is finite as well. 
\end{proof}

For details on derived words of primitive substitutive words, we refer to~\cite{Durand98}.

We remark that, for practical reasons, the set $R_\beta$ can be obtained from the set $R = \{w_\beta\}$ by adding to $R$ iteratively all return words of $0$ which appear in $\psi_\beta(w0)$ for some $w \in R$ until $R$ stabilises.
The final set $R$ is equal to~$R_\beta$.

Now, we apply the theorems in the case of two quadratic examples.

\begin{figure}[ht]
\centering
\begin{tikzpicture}[scale=1.38]
\draw(-4.236,-.1)node[below]{\small$-\beta^3$}--(-4.236,.1) (-3.236,-.1)node[below]{\small$\hspace{.5em}-\beta^3\!+\!\beta^2\!-\!\beta$}--(-3.236,.1) (-2.236,-.1)node[below=2.5ex]{\small$-\beta^3\!+\!\beta^2\!-\!\beta\!+\!1$}--(-2.236,.1) (-1.618,-.1)node[below]{\small$-\beta\vphantom{\beta^2}$}--(-1.618,.1) (-.618,-.1)node[below]{\small$-\beta\!+\!1\vphantom{\beta^2}$}--(-.618,.1) (0,-.1)node[below]{$0\vphantom{\beta^2}$}--(0,.1) (1,-.1)node[below]{$1\vphantom{\beta^2}$}--(1,.1) (2,-.1)node[below]{\small$\beta^2\!-\!\beta\!+\!1$}--(2,.1) (2.618,-.1)node[below=2.5ex]{\small$\beta^2$}--(2.618,.1) (3.618,-.1)node[below]{\small$\beta^4\!-\!\beta^3\!+\!\beta^2\!-\!\beta$}--(3.618,.1) (4.618,-.1)node[below=2.5ex]{\small$\hspace{-1em}\beta^4\!-\!\beta^3\!+\!\beta^2\!-\!\beta\!+\!1$}--(4.618,.1) (5.236,-.1)node[below]{\small$\hspace{1em}\beta^4\!-\!\beta^3\!+\!\beta^2$}--(5.236,.1) (6.236,-.1)node[below=2.5ex]{\small$\hspace{.5em}\beta^4\!-\!\beta\!+\!1$}--(6.236,.1) (6.854,-.1)node[below]{\small$\beta^4$}--(6.854,.1) 
 (-4.236,0)--node[above]{$A$}(-3.236,0)--node[above]{$A$} (-2.236,0)--node[above]{$B$}(-1.618,0)--node[above]{$A$}(-.618,0)--node[above]{$B$}(0,0)--node[above]{$A$}(1,0)--node[above]{$A$}(2,0)--node[above]{$B$}(2.618,0)--node[above]{$A$}(3.618,0)--node[above]{$A$}(4.618,0)--node[above]{$B$}(5.236,0)--node[above]{$A$}(6.236,0)--node[above]{$B$}(6.854,0);
\end{tikzpicture}
\caption{The $(-\beta)$-integers in $[-\beta^3,\beta^4]$, $\beta = (1+\sqrt5)/2$.}
\label{f:intgm}
\end{figure}
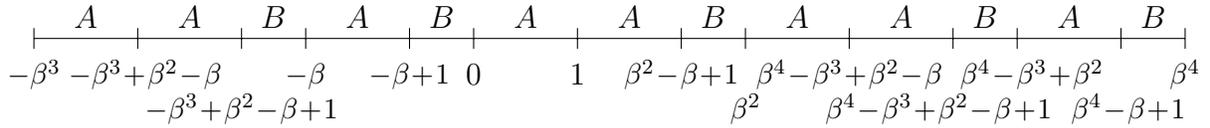

\begin{example} \label{x:gm}
Let $\beta = \frac{1+\sqrt5}{2}$, i.e., $\beta^2 = \beta + 1$, and $t = \frac{-\beta}{\beta+1} = \frac{-1}{\beta}$, then $V_\beta = V'_\beta = \{t, 0\}$.
Since 
\[
J_{\widehat{t}} = (t,0) = \big(t, \tfrac{-1}{\beta^3}\big) \cup \big\{\tfrac{-1}{\beta^3}\big\} \cup \big(\tfrac{-1}{\beta^3}, 0\big)\,, \quad J_{\widehat{0}} = \big(0, \tfrac{1}{\beta^2}\big)\,,
\]
see Figure~\ref{f:Tmbeta} (middle), the anti-morphism $\psi_\beta$ on $A_\beta^*$ is defined by
\[
\psi_\beta: \quad t \mapsto 0 \,, \quad \widehat{t} \mapsto \widehat{0} \, t \, \widehat{t} \,, \quad 0 \mapsto 0 \,, \quad \widehat{0} \mapsto \widehat{t} \,.
\]
Its two-sided fixed point $\cdots u_{-1} u_0 u_1 \cdots$ is
\[
\cdots 
\underbrace{0}_{\psi_\beta(0)} 
\underbrace{\widehat{0}\,t\,\widehat{t}}_{\psi_\beta(\widehat{t}\,)} 
\underbrace{0}_{\psi_\beta(t)} 
\underbrace{\widehat{t}}_{\psi_\beta(\widehat{0})} 
\underbrace{0}_{\psi_\beta(0)} 
\underbrace{\widehat{0}\,t\,\widehat{t}}_{\psi_\beta(\widehat{t}\,)} 
\underbrace{0}_{\psi_\beta(t)} 
\underbrace{\widehat{t}}_{\psi_\beta(\widehat{0})} 
\underbrace{\dot{0}}_{\psi_\beta(\dot{0})} 
\underbrace{\widehat{0}\,t\,\widehat{t}}_{\psi_\beta(\widehat{t}\,)} 
\underbrace{0}_{\psi_\beta(0)} 
\underbrace{\widehat{0}\,t\,\widehat{t}}_{\psi_\beta(\widehat{t}\,)} 
\underbrace{0}_{\psi_\beta(t)} 
\underbrace{\widehat{t}}_{\psi_\beta(\widehat{0})} 
\underbrace{0}_{\psi_\beta(0)} 
\underbrace{\widehat{t}}_{\psi_\beta(\widehat{0})} 
\underbrace{0}_{\psi_\beta(0)}
\cdots \,,
\]
where $\dot{0}$ marks the central letter~$u_0$.
The return word of $0$ starting at~$u_0$ is $w_\beta = 0\, \widehat{0}\, t\, \widehat{t}$.
The image $\psi_\beta(w_\beta 0) = 0 \ \widehat{0} \, t \, \widehat{t} \ 0 \ \widehat{t} \ 0$ contains the return words $w_\beta$ and $0 \, \widehat{t}$. 
Since $\psi_\beta(0 \, \widehat{t}\, 0) = 0 \ \widehat{0} \, t \, \widehat{t} \ 0$, there are no other return words of~$0$, i.e., $R_\beta = \{A, B\}$ with $A = 0 \, \widehat{0} \, t \, \widehat{t}$, $B = 0 \, \widehat{t}$.
Therefore, $\cdots u'_{-2} u'_{-1} u'_1 u'_2 \cdots$ is a two-sided fixed point of the anti-morphism
\[
\varphi_{-\beta}: \quad A \mapsto A B \,,  \quad B \mapsto A \,, 
\]
with
\[
u'_{-13} \cdots u'_{-1} \ u'_1 \cdots u'_{21} = A A B A A B A B A A B A B \ A A B A A B A B A A B A A B A B A A B AB\,.
\]
We have $\lambda(J_{\widehat{0}}) = \frac{1}{\beta^2}$, $\lambda(J_{\widehat{t}}) = \frac{1}{\beta}$, thus $L(A) = 1$, $L(B) = \frac{1}{\beta} = \beta - 1$, and some $(-\beta)$-integers are shown in Figure~\ref{f:intgm}.
Note that $(-\beta)^n$ can also be represented as $(-\beta)^{n+2} + (-\beta)^{n+1}$.
\end{example}

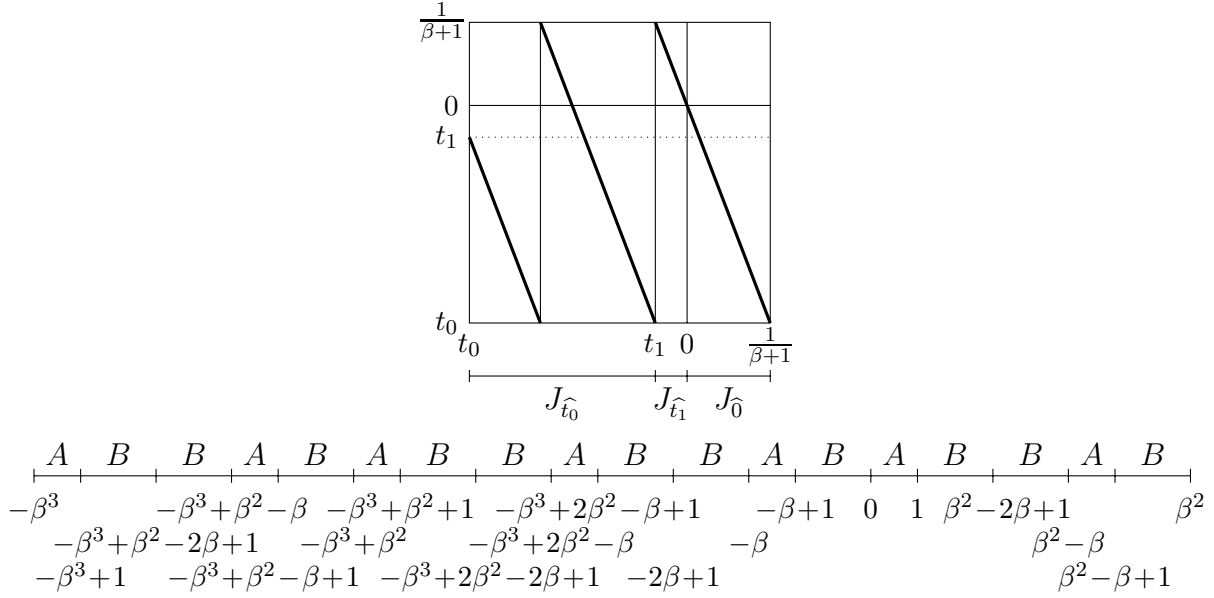
\begin{figure}[ht]
\centering
\begin{tikzpicture}[scale=4]
\draw(-.7236,0)node[left]{\small$0$}--(.2764,0) (0,-.7236)node[below]{\small$0$}--(0,.2764)
(-.7236,-.7236)node[below]{\small$t_0$}node[left]{\small$t_0$}--(.2764,-.7236)node[below=-.5ex]{\small$\frac{1}{\beta+1}$}--(.2764,.2764)--(-.7236,.2764)node[left=-.25em]{\small$\frac{1}{\beta+1}$}--cycle
(-.4875,-.7236)--(-.4875,.2764) (-.1056,-.7236)node[below]{\small$t_1$}--(-.1056,.2764);
\draw[very thick](-.7236,-.1056)--(-.4875,-.7236) (-.4875,.2764)--(-.1056,-.7236) (-.1056,.2764)--(.2764,-.7236);
\draw[dotted](-.7236,-.1056)node[left]{\small$t_1$}--(.2764,-.1056);
\begin{scope}[shift={(0,-.9)}]
\draw(-.7236,-.02)--(-.7236,.02) (-.1056,-.02)--(-.1056,.02) (0,-.02)--(0,.02) (.2764,-.02)--(.2764,.02)
(-.7236,0)--node[below]{$J_{\widehat{t_0}}$}(-.1056,0)--node[below]{$J_{\widehat{t_1}}$}(0,0)--node[below]{$J_{\widehat{0}}$}(.2764,0);
\end{scope}
\end{tikzpicture}

\begin{tikzpicture}[scale=.62]
\draw(-17.944,-.15)node[below]{\small$-\beta^3$}--(-17.944,.15) (-16.944,-.15)node[below=5ex]{\small$-\beta^3\!+\!1$}--(-16.944,.15) (-15.326,-.15)node[below=2.5ex]{\small$-\beta^3\!+\!\beta^2\!-\!2\beta\!+\!1$}--(-15.326,.15) (-13.708,-.15)node[below]{\small$-\beta^3\!+\!\beta^2\!-\!\beta$}--(-13.708,.15)  (-12.708,-.15)node[below=5ex]{\small$\hspace{-1em}-\beta^3\!+\!\beta^2\!-\!\beta\!+\!1$}--(-12.708,.15) (-11.09,-.15)node[below=2.5ex]{\small$-\beta^3\!+\!\beta^2$}--(-11.09,.15) (-10.09,-.15)node[below]{\small$-\beta^3\!+\!\beta^2\!+\!1$}--(-10.09,.15) (-8.472,-.15)node[below=5ex]{\small$\hspace{1em}-\beta^3\!+\!2\beta^2\!-\!2\beta\!+\!1$}--(-8.472,.15) (-6.854,-.15)node[below=2.5ex]{\small$-\beta^3\!+\!2\beta^2\!-\!\beta$}--(-6.854,.15) (-5.854,-.15)node[below]{\small$-\beta^3\!+\!2\beta^2\!-\!\beta\!+\!1$}--(-5.854,.15) (-4.236,-.15)node[below=5ex]{\small$-2\beta\!+\!1\vphantom{\beta^2}$}--(-4.236,.15) (-2.618,-.15)node[below=2.5ex]{\small$-\beta\vphantom{\beta^2}$}--(-2.618,.15) (-1.618,-.15)node[below]{\small$-\beta\!+\!1\vphantom{\beta^2}$}--(-1.618,.15) (0,-.15)node[below]{\small$0\vphantom{\beta^2}$}--(0,.15) (1,-.15)node[below]{\small$1\vphantom{\beta^2}$}--(1,.15) (2.618,-.15)node[below]{\small$\hspace{1em}\beta^2\!-\!2\beta\!+\!1$}--(2.618,.15) (4.236,-.15)node[below=2.5ex]{\small$\beta^2\!-\!\beta$}--(4.236,.15) (5.236,-.15)node[below=5ex]{\small$\beta^2\!-\!\beta\!+\!1$}--(5.236,.15) (6.854,-.15)node[below]{\small$\beta^2$}--(6.854,.15) 
(-17.944,0)--node[above]{$A$}(-16.944,0)--node[above]{$B$}(-15.326,0)--node[above]{$B$}(-13.708,0)--node[above]{$A$}(-12.708,0)--node[above]{$B$}(-11.09,0)--node[above]{$A$}(-10.09,0)--node[above]{$B$}(-8.472,0)--node[above]{$B$}(-6.854,0)--node[above]{$A$}(-5.854,0)--node[above]{$B$}(-4.236,0)--node[above]{$B$}(-2.618,0)--node[above]{$A$}(-1.618,0)--node[above]{$B$}(0,0)--node[above]{$A$}(1,0)--node[above]{$B$}(2.618,0)--node[above]{$B$}(4.236,0)--node[above]{$A$}(5.236,0)--node[above]{$B$}(6.854,0);
\end{tikzpicture}
\caption{The $(-\beta)$-transformation and $\mathbb{Z}_{-\beta} \cap [-\beta^3, \beta^2]$, $\beta = (3+\sqrt5)/2$.}
\label{f:Tgm2}
\end{figure}

\begin{example} \label{x:gm2}
Let $\beta = \frac{3+\sqrt5}{2}$, i.e., $\beta^2 = 3 \beta - 1$, then the $(-\beta)$-transformation is depicted in Figure~\ref{f:Tgm2}, where $t_0 = \frac{-\beta}{\beta+1}$, $t_1 = T_{-\beta}(t_0) = \frac{\beta^2}{\beta+1} - 2 = \frac{-\beta^{-1}}{\beta+1}$, $T_{-\beta}(t_1) = \frac{1}{\beta+1} - 1 = t_0$.
Therefore, $V'_\beta = \{t_0, t_1, 0\}$ and the anti-morphism $\psi_\beta:\, A_\beta^* \to A_\beta^*$ is defined by
\[
\psi_\beta: \quad t_0 \mapsto t_1 \,, \quad \widehat{t_0} \mapsto \widehat{t_0} \, t_1 \widehat{t_1} \, 0 \, \widehat{0} \, t_0 \, \widehat{t_0} \,, \quad t_1 \mapsto t_0 \,, \quad \widehat{t_1} \mapsto \widehat{0} \,, \quad 0 \mapsto 0 \,, \quad \widehat{0} \mapsto \widehat{t_0} \, t_1 \, \widehat{t_1} \,,
\]
which has the two-sided fixed point
\[
\cdots 
\underbrace{0}_{\psi_\beta(0)} 
\underbrace{\widehat{0}}_{\psi_\beta(\widehat{t_1})} 
\underbrace{t_0}_{\psi_\beta(t_1)} 
\underbrace{\widehat{t_0}\,t_1\widehat{t_1}\,0\,\widehat{0}\,t_0\,\widehat{t_0}}_{\psi_\beta(\widehat{t_0})} 
\underbrace{t_1}_{\psi_\beta(t_0)} 
\underbrace{\widehat{t_0}\,t_1\,\widehat{t_1}}_{\psi_\beta(\widehat{0})} 
\underbrace{\dot{0}}_{\psi_\beta(0)} 
\underbrace{\widehat{0}}_{\psi_\beta(\widehat{t_1})} 
\underbrace{t_0}_{\psi_\beta(t_1)} 
\underbrace{\widehat{t_0}\,t_1\widehat{t_1}\,0\,\widehat{0}\,t_0\,\widehat{t_0}}_{\psi_\beta(\widehat{t_0})} 
\cdots \,,
\]
where $\dot{0}$ marks the central letter~$u_0$.
We have $w_\beta = 0 \, \widehat{0} \, t_0 \, \widehat{t_0} \, t_1 \, \widehat{t_1}$ and
\begin{align*}
\psi_\beta: \qquad 
0 \, \widehat{0} \, t_0 \, \widehat{t_0} \, t_1 \, \widehat{t_1} \, 0 & \mapsto 0 \ \widehat{0} \ t_0 \ \widehat{t_0} \, t_1 \, \widehat{t_1} \, 0 \, \widehat{0} \, t_0 \, \widehat{t_0} \ t_1 \ \widehat{t_0} \, t_1 \, \widehat{t_1} \ 0  \,, \\
0 \, \widehat{0} \, t_0 \, \widehat{t_0} \, t_1 \, \widehat{t_0} \, t_1 \, \widehat{t_1} \, 0 & \mapsto 0 \ \widehat{0} \ t_0 \ \widehat{t_0} \, t_1 \, \widehat{t_1} \, 0 \, \widehat{0} \, t_0 \, \widehat{t_0} \ t_0 \ \widehat{t_0} \, t_1 \, \widehat{t_1} \, 0 \, \widehat{0} \, t_0 \, \widehat{t_0} \ t_1 \ \widehat{t_0} \, t_1 \, \widehat{t_1} \ 0  \,, \\
0 \, \widehat{0} \, t_0 \, \widehat{t_0} \, t_0 \, \widehat{t_0} \, t_1 \, \widehat{t_1} \, 0 & \mapsto 0 \ \widehat{0} \ t_0 \ \widehat{t_0} \, t_1 \, \widehat{t_1} \, 0 \, \widehat{0} \, t_0 \, \widehat{t_0} \ t_1 \ \widehat{t_0} \, t_1 \, \widehat{t_1} \, 0 \, \widehat{0} \, t_0 \, \widehat{t_0} \ t_1 \ \widehat{t_0} \, t_1 \, \widehat{t_1} \ 0  \,.
\end{align*}
Note that $0 \, \widehat{0} \, t_0 \, \widehat{t_0} \, t_1 \, \widehat{t_0} \, t_1 \, \widehat{t_1}$ and $0 \, \widehat{0} \, t_0 \, \widehat{t_0} \, t_0 \, \widehat{t_0} \, t_1 \, \widehat{t_1}$ differ only by a letter in~$V'_\beta$, and correspond therefore to intervals of the same length. 
Since the letters $t_0$ and $t_1$ are never mapped to~$0$, we identify these two return words.
This means that $R_\beta = \{A, B\}$ with $A = 0 \, \widehat{0} \, t_0 \, \widehat{t_0} \, t_1 \, \widehat{t_1}$, $B = 0 \, \widehat{0} \, t_0 \, \widehat{t_0} \, \{t_0,t_1\} \, \widehat{t_0} \, t_1 \, \widehat{t_1}$, and 
\[
\cdots u'_{-2} u'_{-1}\ u'_1 u'_2 \cdots = \cdots A B B A B A B B A B B A B \ A B B A B A B B A B B A B \, \cdots
\]
is a two-sided fixed point of the anti-morphism
\[
\varphi_{-\beta}: \quad A \mapsto A B \,,  \quad B \mapsto A B B\,.
\]
We have $L(A) = 1$, $L(B) = \beta - 1 > 1$, and some $(-\beta)$-integers are shown in Figure~\ref{f:Tgm2}.
\end{example}

We remark that it is sufficient to consider the elements of $\,\widehat{\!V'_\beta\!}\,$ when one is only interested in~$\mathbb{Z}_{-\beta}$.
This is made precise in the following proposition.

\begin{proposition} \label{p:hat}
Let $\beta$ and $\cdots u_{-1} u_0 u_1 \cdots$ be as in Theorem~\ref{t:mbeta}, $t = \max\{x \in V_\beta \mid x < 0\}$.
For any $k \in \mathbb{Z}$, $u_{2k} = 0$ is equivalent to $u_{2k-1} = \widehat{t}$ or $u_{2k+1} = \widehat{0}$.

If $0 \not\in V_\beta$ or the size of $V_\beta$ is even, then $u_{2k} = 0$ is equivalent to $u_{2k-1} = \widehat{t}$.
 
If $0 \not\in V_\beta$ or the size of $V_\beta$ is odd, then $u_{2k} = 0$ is equivalent to $u_{2k+1} = \widehat{0}$.
\end{proposition}

\begin{proof}
Let $k \in \mathbb{Z}$ and $m \ge 0$ such that $z_k/\beta^{2m} \in \big(\frac{-\beta}{\beta+1}, \frac{1}{\beta+1}\big)$.
Then we have 
\begin{itemize}
\item
$u_{2k} = 0$ if and only if $T_{-\beta}^{2m}(z_k/\beta^{2m}) = 0$,
\item
$u_{2k-1} = \widehat{t}$ if and only if $\lim_{y\to z_k-} T_{-\beta}^{2m}(y/\beta^{2m}) = 0$,
\item
$u_{2k+1} = \widehat{0}$ if and only if $\lim_{y\to z_k+} T_{-\beta}^{2m}(y/\beta^{2m}) = 0$. 
\end{itemize}
Thus $u_{2k} = 0$, $u_{2k-1} = \widehat{t}$ and $u_{2k+1} = \widehat{0}$ are equivalent when $T_{-\beta}^{2m}$ is continuous at~$z_k/\beta^{2m}$.
Assume from now on that $z_k/\beta^{2m}$ is a discontinuity point of~$T_{-\beta}^{2m}$.
Then $T_{-\beta}^\ell(z_k/\beta^{2m}) = \frac{-\beta}{\beta+1}$ for some $1 \le \ell \le 2m$ and, if $\ell$ is minimal with this property,
\[
\lim_{y\to z_k-} T_{-\beta}^{2\lfloor\ell/2\rfloor+1}(y/\beta^{2m}) = \tfrac{-\beta}{\beta+1} \quad \mbox{and} \quad \lim_{y\to z_k+} T_{-\beta}^{2\lceil\ell/2\rceil}(y/\beta^{2m}) = \tfrac{-\beta}{\beta+1}\,.
\]
Hence, if $0 \not\in V_\beta$, we cannot have $u_{2k} = 0$,  $u_{2k-1} = \widehat{t}$ or $u_{2k+1} = \widehat{0}$ at a discontinuity point, which proves the proposition in this case.
If $0 \in V_\beta$, then $T_{-\beta}^{\#V_\beta-1}\big(\frac{-\beta}{\beta+1}\big) = 0$, thus
\begin{itemize}
\item
$T_{-\beta}^j(z_k/\beta^{2m}) = 0$ if and only if $j \ge \ell + \#V_\beta - 1$,
\item
$\lim_{y\to z_k-} T_{-\beta}^j(y/\beta^{2m}) = 0$ if and only if $j \ge 2\lfloor\ell/2\rfloor + \#V_\beta$,
\item
$\lim_{y\to z_k+} T_{-\beta}^j(y/\beta^{2m}) = 0$ if and only if $j \ge 2\lceil\ell/2\rceil + \#V_\beta - 1$.
\end{itemize}
Since $2\lfloor\ell/2\rfloor \ge \ell-1$ and $2\lceil\ell/2\rceil \ge \ell$, we obtain $u_{2k} = 0$ whenever $u_{2k-1} = \widehat{t}$ or $u_{2k+1} = \widehat{0}$.
If $\#V_\beta$ is even, then $u_{2k} = 0$ implies that $u_{2k-1} = \widehat{t}$ since $2m \ge \ell + \#V_\beta - 1$ implies that $2m \ge 2\lfloor\ell/2\rfloor + \#V_\beta$.
If $\#V_\beta$ is odd, then $u_{2k} = 0$ implies that $u_{2k+1} = \widehat{0}$ since $2m \ge \ell + \#V_\beta - 1$ implies that $2m \ge 2\lceil\ell/2\rceil + \#V_\beta - 1$. 
This proves the proposition.
\end{proof}

By Proposition~\ref{p:hat}, it suffices to consider the anti-morphism $\widehat{\psi}_\beta:\, \,\widehat{\!V'_\beta\!}\,^* \to \,\widehat{\!V'_\beta\!}\,^*$ defined by
\[
\widehat{\psi}_\beta(\widehat{x}) = \,\widehat{\!x_m\!}\, \cdots \widehat{x_1} \, \widehat{x_0} \quad \mbox{when} \quad \psi_\beta(\widehat{x}) = \,\widehat{\!x_m\!}\, \, T_{-\beta}(y_m) \cdots \widehat{x_1} \, T_{-\beta}(y_1) \, \widehat{x_0} \quad (x \in V'_\beta).
\]
Then $\Delta_{-\beta}$ is given by the set $\widehat{R}_\beta$ which consists of the return words of $\widehat{0}$ when $0 \not\in V_\beta$ or the size of $V_\beta$ is odd.
When $0 \in V_\beta$ and the size of $V_\beta$ is even, as in Example~\ref{x:gm}, then $\widehat{R}_\beta$ consists of the words $w\, \widehat{t}$ such that $\widehat{t}\, w$ is a return word of $\widehat{t}$. 

\begin{figure}[ht]
\centering
\begin{tikzpicture}[scale=6]
\draw(-.688,0)node[left]{$0$}--(.312,0) (0,-.688)node[below]{$0$}--(0,.312)
(-.688,-.688)node[below]{$t_0$}node[left]{$t_0$}--(.312,-.688)node[below=-.5ex]{$\frac{1}{\beta+1}$}--(.312,.312)--(-.688,.312)node[left=-.25em]{$\frac{1}{\beta+1}$}--cycle
(-.5948,-.688)--(-.5948,.312) (-.1414,-.688)node[below]{$t_3$}--(-.1414,.312);
\draw[very thick](-.688,-.4825)--(-.5948,-.688) (-.5948,.312)--(-.1414,-.688) (-.1414,.312)--(.312,-.688);
\draw[dotted](-.688,-.4825)node[left]{$t_1$}--(.312,-.4825) (-.4825,-.688)node[below]{$t_1$}--(-.4825,.312) (-.688,.0641)node[left]{$t_2$}--(.312,.0641) (.0641,-.688)node[below]{$t_2$}--(.0641,.312) (-.688,-.1414)node[left]{$t_3$}--(.312,-.1414); 
\begin{scope}[shift={(0,-.81)}]
\draw(-.688,-.015)--(-.688,.015) (-.4825,-.015)--(-.4825,.015) (-.1414,-.015)--(-.1414,.015) (0,-.015)--(0,.015) (.0641,-.015)--(.0641,.015) (.312,-.015)--(.312,.015)
(-.688,0)--node[below]{$J_{\widehat{t_0}}$}(-.4825,0)--node[below]{$J_{\widehat{t_1}}$}(-.1414,0)--node[below]{$J_{\widehat{t_3}}$}(0,0)--node[below]{$J_{\widehat{0}}$}(.0641,0)--node[below]{$J_{\widehat{t_2}}$}(.312,0);
\end{scope}
\end{tikzpicture}

\begin{tikzpicture}[scale=.445]
\draw(-10.73,-.2)node[below]{\small$-\beta^3$}--(-10.73,.2) (-9.73,-.2)node[below=7.5ex]{\small$-\beta^3\!+\!1$}--(-9.73,.2) (-8.07,-.2)node[below=5ex]{\small$-\beta^3\!+\!\beta^2\!-\!\beta$}--(-8.07,.2) 
(-7.07,-.2)node[below=2.5ex]{\small$-\beta^3\!+\!\beta^2\!-\!\beta\!+\!1$}--(-7.07,.2) 
(-5.86,-.2)node[below]{\small$-\beta^3\!+\!\beta^2$}--(-5.86,.2) 
(-4.86,-.2)node[below=7.5ex]{\small$-\beta^3\!+\!\beta^2\!+\!1$}--(-4.86,.2) 
(-2.21,-.2)node[below=5ex]{\small$-\beta$}--(-2.21,.2) 
(-1.21,-.2)node[below=2.5ex]{\small$-\beta\!+\!1$}--(-1.21,.2) 
(0,-.2)node[below]{\small$\vphantom{\beta^2}0$}--(0,.2) 
(1,-.2)node[below]{\small$\vphantom{\beta^2}1$}--(1,.2) 
(2.66,-.2)node[below=2.5ex]{\small$\beta^2\!-\!\beta$}--(2.66,.2) 
(3.66,-.2)node[below=5ex]{\small$\beta^2\!-\!\beta\!+\!1$}--(3.66,.2) 
(4.86,-.2)node[below=7.5ex]{\small$\beta^2$}--(4.86,.2) 
(5.86,-.2)node[below]{\small$\beta^2\!+\!1$}--(5.86,.2) 
(8.07,-.2)node[below=2.5ex]{\small$\beta^4\!-\!2\beta^3\!+\!\beta^2\!+\!1$}--(8.07,.2) 
(10.73,-.2)node[below=5ex]{\small$\hspace{-1em}\beta^4\!-\!\beta^3\!-\!\beta$}--(10.73,.2) 
(11.73,-.2)node[below=7.5ex]{\small$\hspace{-2em}\beta^4\!-\!\beta^3\!-\!\beta\!+\!1$}--(11.73,.2) 
(12.93,-.2)node[below]{\small$\beta^4\!-\!\beta^3$}--(12.93,.2) 
(13.93,-.2)node[below=2.5ex]{\small$\hspace{-.5em}\beta^4\!-\!\beta^3\!+\!1$}--(13.93,.2) 
(15.59,-.2)node[below=5ex]{\small$\beta^4\!-\!\beta^3\!+\!\beta^2\!-\!\beta$}--(15.59,.2) 
(16.59,-.2)node[below=7.5ex]{\small$\hspace{1em}\beta^4\!-\!\beta^3\!+\!\beta^2\!-\!\beta\!+\!1$}--(16.59,.2) 
(17.8,-.2)node[below]{\small$\beta^4\!-\!\beta^3\!+\!\beta^2$}--(17.8,.2) 
(18.8,-.2)node[below=2.5ex]{\small$\beta^4\!-\!\beta^3\!+\!\beta^2\!+\!1$}--(18.8,.2) 
(21.46,-.2)node[below=5ex]{\small$\beta^4\!-\!\beta$}--(21.46,.2) 
(22.46,-.2)node[below=7.5ex]{\small$\beta^4\!-\!\beta\!+\!1$}--(22.46,.2) 
(23.66,-.2)node[below]{\small$\beta^4$}--(23.66,.2) 
(-10.73,0)--node[above]{$A$}(-9.73,0)--node[above]{$C$}(-8.07,0)--node[above]{$A$}(-7.07,0)--node[above]{$B$}(-5.86,0)--node[above]{$A$}(-4.86,0)--node[above]{$D$}(-2.21,0)--node[above]{$A$}(-1.21,0)--node[above]{$B$}(0,0)--node[above]{$A$}(1,0)--node[above]{$C$}(2.66,0)--node[above]{$A$}(3.66,0)--node[above]{$B$}(4.86,0)--node[above]{$A$}(5.86,0)--node[above]{$E$}(8.07,0)--node[above]{$D$}(10.73,0)--node[above]{$A$}(11.73,0)--node[above]{$B$}(12.93,0)--node[above]{$A$}(13.93,0)--node[above]{$C$}(15.59,0)--node[above]{$A$}(16.59,0)--node[above]{$B$}(17.8,0)--node[above]{$A$}(18.8,0)--node[above]{$D$}(21.46,0)--node[above]{$A$}(22.46,0)--node[above]{$B$}(23.66,0);
\end{tikzpicture}
\caption{The $(-\beta)$-transformation and $\mathbb{Z}_{-\beta} \cap [-\beta^3,\beta^4]$ from Example~\ref{x:complex}.}
\label{f:Tcomplex}
\end{figure}
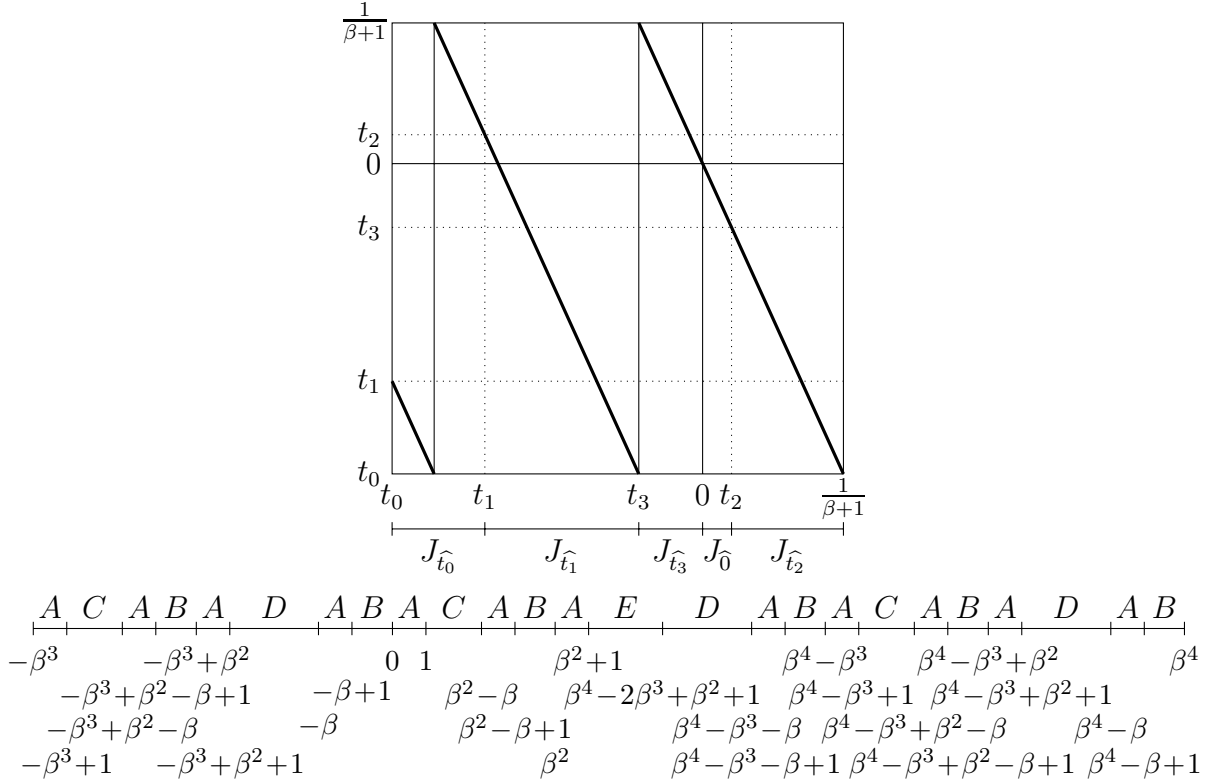

\begin{example} \label{x:complex}
Let $\beta > 1$ with $\beta^3 = 2 \beta^2 + 1$, i.e., $\beta \approx 2.206$, and let $t_n = T_{-\beta}^n\big(\frac{-\beta}{\beta+1}\big)$ for $n \ge 0$.
Then we have 
\[
t_0 = \tfrac{-\beta}{\beta+1}, \quad t_1 = \tfrac{\beta^2}{\beta+1} - 2 = \tfrac{\beta^{-1}-2}{\beta+1}, \quad t_2 = \tfrac{2\beta-1}{\beta+1} - 1 = \tfrac{\beta^{-2}}{\beta+1},  \quad t_3 = \tfrac{-\beta^{-1}}{\beta+1}, \quad t_4 = \tfrac{1}{\beta+1} - 1 = t_0,
\]
see Figure~\ref{f:Tcomplex}.
The anti-morphism $\widehat{\psi}_\beta:\, \,\widehat{\!V'_\beta\!}\,^* \to \,\widehat{\!V'_\beta\!}\,^*$ is therefore defined by
\[
\widehat{\psi}_\beta: \quad 
\widehat{t_0} \mapsto \widehat{t_2} \, \widehat{t_0} \,, \quad
\widehat{t_1} \mapsto \widehat{t_0} \, \widehat{t_1} \, \widehat{t_3} \, \widehat{0} \,, \quad
\widehat{t_3} \mapsto \widehat{0} \, \widehat{t_2} \,, \quad
\widehat{0} \mapsto \widehat{t_3} \,, \quad
\widehat{t_2} \mapsto \widehat{t_0} \, \widehat{t_1} \,.
\]
Since $0 \not\in V_\beta$, we consider return words of~$\widehat{0}$ in the $\widehat{\psi}_\beta$-images of $\widehat{w}_\beta = \widehat{0} \, \widehat{t_2} \, \widehat{t_0} \, \widehat{t_1} \, \widehat{t_3}$:
\begin{align*}
\widehat{0} \, \widehat{t_2} \, \widehat{t_0} \, \widehat{t_1} \, \widehat{t_3} & \mapsto \widehat{0} \, \widehat{t_2} \, \widehat{t_0} \, \widehat{t_1} \, \widehat{t_3} \ \widehat{0} \, \widehat{t_2} \, \widehat{t_0} \, \widehat{t_0} \, \widehat{t_1} \, \widehat{t_3} \,, \\
\widehat{0} \, \widehat{t_2} \, \widehat{t_0} \, \widehat{t_0} \, \widehat{t_1} \, \widehat{t_3} & \mapsto \widehat{0} \, \widehat{t_2} \, \widehat{t_0} \, \widehat{t_1} \, \widehat{t_3} \ \widehat{0} \, \widehat{t_2} \, \widehat{t_0} \, \widehat{t_2} \, \widehat{t_0} \, \widehat{t_0} \, \widehat{t_1} \, \widehat{t_3} \,, \\
\widehat{0} \, \widehat{t_2} \, \widehat{t_0} \, \widehat{t_2} \, \widehat{t_0} \, \widehat{t_0} \, \widehat{t_1} \, \widehat{t_3} & \mapsto \widehat{0} \, \widehat{t_2} \, \widehat{t_0} \, \widehat{t_1} \, \widehat{t_3} \ \widehat{0} \, \widehat{t_2} \, \widehat{t_0} \, \widehat{t_2} \, \widehat{t_0} \, \widehat{t_0} \, \widehat{t_1} \, \widehat{t_2} \, \widehat{t_0} \, \widehat{t_0} \, \widehat{t_1} \, \widehat{t_3} \,, \\
\widehat{0} \, \widehat{t_2} \, \widehat{t_0} \, \widehat{t_2} \, \widehat{t_0} \, \widehat{t_0} \, \widehat{t_1} \, \widehat{t_2} \, \widehat{t_0} \, \widehat{t_0} \, \widehat{t_1} \, \widehat{t_3} & \mapsto \widehat{0} \, \widehat{t_2} \, \widehat{t_0} \, \widehat{t_1} \, \widehat{t_3} \ \widehat{0} \, \widehat{t_2} \, \widehat{t_0} \, \widehat{t_2} \, \widehat{t_0} \, \widehat{t_0} \, \widehat{t_1} \, \widehat{t_0} \, \widehat{t_1} \, \widehat{t_3} \ \widehat{0} \, \widehat{t_2} \, \widehat{t_0} \, \widehat{t_2} \, \widehat{t_0} \, \widehat{t_0} \, \widehat{t_1} \, \widehat{t_2} \, \widehat{t_0} \, \widehat{t_0} \, \widehat{t_1} \, \widehat{t_3} \,, \\
\widehat{0} \, \widehat{t_2} \, \widehat{t_0} \, \widehat{t_2} \, \widehat{t_0} \, \widehat{t_0} \, \widehat{t_1} \, \widehat{t_0} \, \widehat{t_1} \, \widehat{t_3} & \mapsto \widehat{0} \, \widehat{t_2} \, \widehat{t_0} \, \widehat{t_1} \, \widehat{t_3} \ \widehat{0} \, \widehat{t_2} \, \widehat{t_0} \, \widehat{t_0} \, \widehat{t_1} \, \widehat{t_3}\ \widehat{0} \, \widehat{t_2} \, \widehat{t_0} \, \widehat{t_2} \, \widehat{t_0} \, \widehat{t_0} \, \widehat{t_1} \, \widehat{t_2} \, \widehat{t_0} \, \widehat{t_0} \, \widehat{t_1} \, \widehat{t_3}\,.
\end{align*}
Hence $\widehat{R}_\beta = \{A, B, C, D, E\}$ with $A = \widehat{0} \, \widehat{t_2} \, \widehat{t_0} \, \widehat{t_1} \, \widehat{t_3}$, $B = \widehat{0} \, \widehat{t_2} \, \widehat{t_0} \, \widehat{t_0} \, \widehat{t_1} \, \widehat{t_3}$, $C = \widehat{0} \, \widehat{t_2} \, \widehat{t_0} \, \widehat{t_2} \, \widehat{t_0} \, \widehat{t_0} \, \widehat{t_1} \, \widehat{t_3}$, $D = \widehat{0} \, \widehat{t_2} \, \widehat{t_0} \, \widehat{t_2} \, \widehat{t_0} \, \widehat{t_0} \, \widehat{t_1} \, \widehat{t_2} \, \widehat{t_0} \, \widehat{t_0} \, \widehat{t_1} \, \widehat{t_3}$, $E = \widehat{0} \, \widehat{t_2} \, \widehat{t_0} \, \widehat{t_2} \, \widehat{t_0} \, \widehat{t_0} \, \widehat{t_1} \, \widehat{t_0} \, \widehat{t_1} \, \widehat{t_3}$, and $\mathbb{Z}_{-\beta}$ is described by the anti-morphism $\widehat{\varphi}_{-\beta}:\, \widehat{R}_\beta^* \to \widehat{R}_\beta^*$ given by
\[
\widehat{\varphi}_{-\beta}: \quad A \mapsto A B\,, \quad B \mapsto A C\,, \quad C \mapsto A D\,, \quad D \mapsto A E D\,, \quad E \mapsto A B D\,.
\]
The $(-\beta)$-integers in $[-\beta^3, \beta^4]$ are represented in Figure~\ref{f:Tcomplex}, and we have 
\[
L(A) = 1, \quad L(B) = \beta - 1, \quad L(C) = \beta^2 - \beta - 1, \quad L(D) = \beta^2 - \beta \approx 2.659, \quad L(E) = \beta. 
\]
Note that $L(D) > \beta > 2$.
Moreover, the cardinality of $\Delta_{-\beta}$ is larger than that of~$V_\beta$, which in turn is larger than the algebraic degree $d$ of~$\beta$ ($\#\Delta_{-\beta} = 5$, $\#V_\beta = 4$, $d = 3$). 
\end{example}

\begin{example} \label{x:complex2}
Let $\beta > 1$ with $\beta^6 = 3 \beta^5 + 2 \beta^4 + 2 \beta^3 + \beta^2 - 2 \beta - 1$, i.e., $\beta \approx 3.695$, then the $(-\beta)$-transformation is depicted in Figure~\ref{f:Tcomplex2}, where $t_n = T_{-\beta}^n\big(\frac{-\beta}{\beta+1}\big)$.
We have $t_5 = \frac{-1}{\beta+1} = t_6$.
The anti-morphism $\widehat{\psi}_\beta:\, \,\widehat{\!V'_\beta\!}\,^* \to \,\widehat{\!V'_\beta\!}\,^*$ is therefore given by
\begin{align*}
\widehat{\psi}_\beta: \quad 
\widehat{t_0} & \mapsto \widehat{t_3} \, \widehat{t_5} \,, &
\widehat{t_2} & \mapsto \widehat{t_4} \, \widehat{t_0} \, \widehat{t_2} \,, & 
\widehat{t_3} & \mapsto \widehat{t_5} \, \widehat{t_1} \, \widehat{0} \, \widehat{t_4} \, \widehat{t_0} \, \widehat{t_2} \, \widehat{t_3} \, \widehat{t_5} \, \widehat{t_1} \, \widehat{0} \,, \hspace{-6em} \\
 \widehat{t_5} & \mapsto \widehat{t_2} \, \widehat{t_3} \,, &
\widehat{t_1} & \mapsto \widehat{0} \, \widehat{t_4} \, \widehat{t_0} \,, & 
\widehat{0} & \mapsto \widehat{t_5} \, \widehat{t_1} \,, &
\widehat{t_4} & \mapsto \widehat{t_0} \, \widehat{t_2} \, \widehat{t_3} \,.
\end{align*}
In order to deal with shorter words, we group the letters forming the words 
\[
a = \widehat{0} \, \widehat{t_4},\, 
b = \widehat{t_0} \, \widehat{t_2} \, \widehat{t_3} \, \widehat{t_5} \, \widehat{t_1},\, 
c = \widehat{t_0} \, \widehat{t_2} \, \widehat{t_3} \, \widehat{t_5},\, 
d = \widehat{t_2} \, \widehat{t_3} \, \widehat{t_5} \, \widehat{t_1},\, 
e = \widehat{t_0} \, \widehat{t_2},\,
f\! = \widehat{t_4},\, 
g = \widehat{t_0} \, \widehat{t_2} \, \widehat{t_3},\, 
h = \widehat{t_5} \, \widehat{t_1}, 
\]
which correspond to the intervals $J_a = \big(0,\tfrac{1}{\beta+1}\big)$, $J_b = (t_0,0)$, $J_c = (t_0,t_1)$, $J_d = (t_2,0)$, $J_e = (t_0,t_3)$, $J_f = \big(t_4,\tfrac{1}{\beta+1}\big)$, $J_g = (t_0,t_5)$, $J_h= (t_5,0)$, occurring in iterated images of~$J_a$.
The anti-morphism $\widehat{\psi}_\beta$ acts on these words by
\begin{align*}
\widehat{\psi}_\beta: \quad 
a & \mapsto b \,, &
b & \mapsto a b a b a c \,, & 
c & \mapsto d a b a c \,, &
d & \mapsto a b a b a e \,, \\
e & \mapsto f c \,, &
f & \mapsto g \,, & 
g & \mapsto h a b a c \,, &
h & \mapsto a g \,.
\end{align*}
Since $\widehat{0}$ only occurs at the beginning of~$a$, the return words of $\widehat{0}$ with their $\widehat{\psi}_\beta$-images are 
\begin{align*}
a b & \mapsto a b \ a b \ a c b \,, &
a e d & \mapsto a b \ a b \ a e f c b \,, \\
a c b & \mapsto a b \ a b \ a c d \ a b \ a c b \,, &
a e f c b & \mapsto a b \ a b \ a c d \ a b \ a c g f c b \,, \\
a c d & \mapsto a b \ a b \ a e d \ a b \ a c b \,, &
a c g f c b & \mapsto a b \ a b \ a c d \ a b \underbrace{a c g h}_{\displaystyle=\!a c b} a b \ a c d \ a b \ a c b \,.
\end{align*}
Therefore, $\mathbb{Z}_{-\beta}$ is described by the anti-morphism $\widehat{\varphi}_{-\beta}:\, \widehat{R}_\beta^* \to \widehat{R}_\beta^*$ which is defined by
\begin{align*}
\widehat{\varphi}_{-\beta}: \quad 
A & \mapsto A A B  \,, & L(A) & = 1 \,, \\
B & \mapsto A A C A B \,, & L(B) & = \beta - 2 \approx 1.695 \,, \\
C & \mapsto A A D A B \,, & L(C) & = \beta^2 - 3 \beta - 1 \approx 1.569 \,, \\
D & \mapsto A A E \,, & L(D) & = \beta^3 - 3 \beta^2 - 2 \beta - 1 \approx 1.104 \,, \\
E & \mapsto A A C A F \,, & L(E) & = \beta^4 - 3 \beta^3 - 2 \beta^2 - \beta - 2 \approx 2.081 \,, \\
F & \mapsto A A C A B A C A B \,, & L(F) & = \beta^5 - 3 \beta^4 - 2 \beta^3 - 2 \beta^2 + \beta -2 \approx 3.12 \,.
\end{align*}
\end{example}

\begin{figure}[ht]
\centering
\begin{tikzpicture}[scale=6]
\draw(-.787,0)node[left]{\small$0$}--(.213,0) (0,-.787)node[below]{\small$0$}--(0,.213)
(-.787,-.787)node[below]{\small$t_0$}node[left]{\small$t_0$}--(.213,-.787)node[below=-.5ex]{\small$\frac{1}{\beta+1}$}--(.213,.213)--(-.787,.213)node[left=-.25em]{\small$\frac{1}{\beta+1}$}--cycle
(-.5989,-.787)--(-.5989,.213) (-.3282,-.787)--(-.3282,.213) (-.0576,-.787)--(-.0576,.213);
\draw[very thick](-.787,-.0917)--(-.5989,-.787) (-.5989,.213)--(-.3282,-.787) (-.3282,.213)--(-.0576,-.787) (-.0576,.213)--(.213,-.787);
\draw[dotted](-.787,-.0917)node[left]{\small$t_1$}--(.213,-.0917) (-.0917,-.787)node[below]{\small$t_1$}--(-.0917,.213) (-.787,-.6615)node[left]{\small$t_2$}--(.213,-.6615) (-.6615,-.787)node[below]{\small$t_2$}--(-.6615,.213) (-.787,-.5567)node[left]{\small$t_3$}--(.213,-.5567) (-.5567,-.787)node[below]{\small$t_3$}--(-.5567,.213)  (-.787,.0576)node[left]{\small$t_4$}--(.213,.0576) (.0576,-.787)node[below]{\small$t_4$}--(.0576,.213)  (-.787,-.213)node[left]{\small$t_5$}--(.213,-.213) (-.213,-.787)node[below]{\small$t_5$}--(-.213,.213); 
\begin{scope}[shift={(0,-.9)}]
\draw(-.787,-.015)--(-.787,.015) (-.6615,-.015)--(-.6615,.015) (-.5567,-.015)--(-.5567,.015) (-.213,-.015)--(-.213,.015) (-.0917,-.015)--(-.0917,.015) (0,-.015)--(0,.015) (.0576,-.015)--(.0576,.015) (.213,-.015)--(.213,.015)
(-.787,0)--node[below]{$J_{\widehat{t_0}}$}(-.6615,0)--node[below]{$J_{\widehat{t_2}}$}(-.5567,0)--node[below]{$J_{\widehat{t_3}}$}(-.213,0)--node[below]{$J_{\widehat{t_5}}$}(-.0917,0)--node[below]{$J_{\widehat{t_1}}$}(0,0)--node[below]{$J_{\widehat{0}}$}(.0576,0)--node[below]{$J_{\widehat{t_4}}$}(.213,0);
\end{scope}
\end{tikzpicture}

\begin{tikzpicture}[scale=.89]
\draw(-3.695,-.15)node[below]{\small$-\beta\vphantom{\beta^2}$}--(-3.695,.15) (-2.695,-.15)node[below=2.5ex]{\small$-\beta\!+\!1\vphantom{\beta^2}$}--(-2.695,.15) (-1.695,-.15)node[below]{\small$-\beta\!+\!2\vphantom{\beta^2}$}--(-1.695,.15) (0,-.15)node[below]{\small$0\vphantom{\beta^2}$}--(0,.15) (1,-.15)node[below]{\small$1\vphantom{\beta^2}$}--(1,.15) (2,-.15)node[below]{\small$2\vphantom{\beta^2}$}--(2,.15) (3.569,-.15)node[below]{\small$\beta^2\!-\!3\beta\!+\!1$}--(3.569,.15) (4.569,-.15)node[below=2.5ex]{\small$\beta^2\!-\!3\beta\!+\!2$}--(4.569,.15) (6.265,-.15)node[below]{\small$\beta^2\!-\!2\beta$}--(6.265,.15) (7.265,-.15)node[below=2.5ex]{\small$\beta^2\!-\!2\beta\!+\!1$}--(7.265,.15) (8.265,-.15)node[below]{\small$\beta^2\!-\!2\beta\!+\!2$}--(8.265,.15) (9.96,-.15)node[below=2.5ex]{\small$\beta^2\!-\!\beta$}--(9.96,.15) (10.96,-.15)node[below]{\small$\beta^2\!-\!\beta\!+\!1$}--(10.96,.15) 
(11.96,-.15)node[below=2.5ex]{\small$\beta^2\!-\!\beta\!+\!2$}--(11.96,.15) (13.66,-.15)node[below]{\small$\beta^2$}--(13.66,.15) 
(-3.695,0)--node[above]{$A$}(-2.695,0)--node[above]{$A$}(-1.695,0)--node[above]{$B$}(0,0)--node[above]{$A$}(1,0)--node[above]{$A$}(2,0)--node[above]{$C$}(3.569,0)--node[above]{$A$}(4.569,0)--node[above]{$B$}(6.265,0)--node[above]{$A$}(7.265,0)--node[above]{$A$}(8.265,0)--node[above]{$B$}(9.96,0)--node[above]{$A$}(10.96,0)--node[above]{$A$}(11.96,0)--node[above]{$B$}(13.66,0);
\end{tikzpicture}
\caption{The $(-\beta)$-transformation and $\mathbb{Z}_{-\beta} \cap [-\beta,\beta^2]$ from Example~\ref{x:complex2}.}
\label{f:Tcomplex2}
\end{figure}

\section{Conclusions}
With every Yrrap number $\beta \ge (1+\sqrt5)/2$, we have associated an anti-morphism $\varphi_{-\beta}$ on a finite alphabet.
The distances between consecutive $(-\beta)$-integers are described by a fixed point of~$\varphi_{-\beta}$.
In~\cite{Ambroz-Dombek-Masakova-Pelantova}, the anti-morphism is described explicitely for each $\beta > 1$ such that $T_{-\beta}^n\big(\frac{-\beta}{\beta+1}\big) \le 0$ and $T_{-\beta}^{2n-1}\big(\frac{-\beta}{\beta+1}\big) \ge \frac{1-\lfloor\beta\rfloor}{\beta}$ for all $n \ge 1$.
Examples~\ref{x:complex} and~\ref{x:complex2} show that the situation can be quite complicated when this condition is not fulfilled. 
Although $\varphi_{-\beta}$ can be obtained by a simple algorithm, it seems to be difficult to find a priori bounds for the number of different distances between consecutive $(-\beta)$-integers or for their maximal~value. 
Only the case of quadratic Pisot numbers~$\beta$ is completely solved; here, we know from \cite{Ito-Sadahiro09,Ambroz-Dombek-Masakova-Pelantova} that $\#V_\beta = \#\Delta_{-\beta} = 2$.

Recall that the maximal distance between consecutive $\beta$-integers is~$1$, and the number of different distances is equal to the cardinality of the set $\{T_\beta^n(1^-) \mid n \ge 0\}$.
Example~\ref{x:complex} shows that the $(-\beta)$-integers do not satisfy similar properties.
By generalising Example~\ref{x:complex2} to $\beta > 1$ with $\beta^6 = (m+1) \beta^5 + m \beta^4 + m \beta^3 + \beta^2 - m \beta - 1$, $m \ge 2$, one sees that the maximal distance can be arbitrarily close to~$4$ for algebraic integers of degree~$6$ and $\#V_\beta = 6$.

In a forthcoming paper, we associate anti-morphisms~$\varphi_{-\beta}$ on infinite alphabets with non-Yrrap numbers~$\beta$, by considering the intervals occurring in the iterated $T_{-\beta}$-images of $\big(0, \frac{1}{\beta+1}\big)$, cf.\ Example~\ref{x:complex2}, and we show that the distances between consecutive $(-\beta)$-integers can be unbounded, e.g.\ for $\beta > 1$ satisfying $\frac{-\beta}{\beta+1} = \sum_{k=1}^\infty a_k (-\beta)^{-k}$ where $a_1 a_2 \cdots = 31232\, 1\, 2\, 31232\, 2\, \cdots$ is a fixed point of the morphism $3 \mapsto 31232,\ 2 \mapsto 2,\ 1 \mapsto 1$.
For Yrrap numbers~$\beta$, this implies that there is no bound for the distance between consecutive $(-\beta)$-integers which is independent of~$\beta$.
However, large distances occur probably only far away from~$0$ and when $\#V_\beta$ is large, and it would be interesting to quantify these relations. 

Another topic that is worth investigating is the structure of the sets $S_{-\beta}(x)$ for $x \ne 0$, and of the corresponding tilings when $\beta$ is a Pisot unit.
A~related question is whether $\mathbb{Z}_{-\beta}$ can be given by a cut and project scheme, cf.~\cite{Berthe-Siegel05,VergerGaugry-Gazeau04}.

\bibliographystyle{amsplain}
\bibliography{betaintegers}

\providecommand{\bysame}{\leavevmode\hbox to3em{\hrulefill}\thinspace}
\providecommand{\MR}{\relax\ifhmode\unskip\space\fi MR }
\providecommand{\MRhref}[2]{%
  \href{http://www.ams.org/mathscinet-getitem?mr=#1}{#2}
}
\providecommand{\href}[2]{#2}
\begin{thebibliography}{10}

\bibitem{Ambroz-Dombek-Masakova-Pelantova}
P.~Ambro{\v{z}}, D.~Dombek, Z.~Mas{\'a}kov{\'a}, and E.~Pelantov{\'a},
  \emph{Numbers with integer expansion in the numeration system with negative
  base}, arXiv:0912.4597v3.

\bibitem{Balkova-Gazeau-Pelantova08}
L.~Balkov{\'a}, J.-P. Gazeau, and E.~Pelantov{\'a}, \emph{Asymptotic behavior
  of beta-integers}, Lett. Math. Phys. \textbf{84} (2008), no.~2-3, 179--198.

\bibitem{Balkova-Pelantova-Steiner08}
L.~Balkov{\'a}, E.~Pelantov{\'a}, and W.~Steiner, \emph{Sequences with constant
  number of return words}, Monatsh. Math. \textbf{155} (2008), no.~3-4,
  251--263.

\bibitem{Bernat-Masakova-Pelantova07}
J.~Bernat, Z.~Mas{\'a}kov{\'a}, and E.~Pelantov{\'a}, \emph{On a class of
  infinite words with affine factor complexity}, Theoret. Comput. Sci.
  \textbf{389} (2007), no.~1-2, 12--25.

\bibitem{Berthe-Siegel05}
V.~Berth{\'e} and A.~Siegel, \emph{Tilings associated with beta-numeration and
  substitutions}, Integers \textbf{5} (2005), no.~3, A2, 46 pp. (electronic).

\bibitem{Burdik-Frougny-Gazeau-Krejcar98}
{\v{C}}.~Burd{\'{\i}}k, C.~Frougny, J.~P. Gazeau, and R.~Krejcar,
  \emph{Beta-integers as natural counting systems for quasicrystals}, J. Phys.
  A \textbf{31} (1998), no.~30, 6449--6472.

\bibitem{Durand98}
F.~Durand, \emph{A characterization of substitutive sequences using return
  words}, Discrete Math. \textbf{179} (1998), no.~1-3, 89--101.

\bibitem{Enomoto08}
F.~Enomoto, \emph{{$AH$}-substitution and {M}arkov partition of a group
  automorphism on {$T^d$}}, Tokyo J. Math. \textbf{31} (2008), no.~2, 375--398.

\bibitem{Fabre95}
S.~Fabre, \emph{Substitutions et {$\beta$}-syst\`emes de num\'eration},
  Theoret. Comput. Sci. \textbf{137} (1995), no.~2, 219--236.

\bibitem{Frougny-Lai09}
C.~Frougny and A.~C. Lai, \emph{On negative bases}, Proceedings of {DLT} 09,
  Lecture Notes in Comput. Sci., vol. 5583, Springer, Berlin, 2009,
  pp.~252--263.

\bibitem{Frougny-Masakova-Pelantova04}
C.~Frougny, Z.~Mas{\'a}kov{\'a}, and E.~Pelantov{\'a}, \emph{Complexity of
  infinite words associated with beta-expansions}, Theor. Inform. Appl.
  \textbf{38} (2004), no.~2, 163--185, Corrigendum: {T}heor. {I}nform. {A}ppl.
  {\bf 38} (2004), no. 3, 269--271.

\bibitem{VergerGaugry-Gazeau04}
J.-P. Gazeau and J.-L. Verger-Gaugry, \emph{Geometric study of the
  beta-integers for a {P}erron number and mathematical quasicrystals}, J.
  Th\'eor. Nombres Bordeaux \textbf{16} (2004), no.~1, 125--149.

\bibitem{Gora07}
P.~G{\'o}ra, \emph{Invariant densities for generalized {$\beta$}-maps}, Ergodic
  Theory Dynam. Systems \textbf{27} (2007), no.~5, 1583--1598.

\bibitem{Ito-Sadahiro09}
S.~Ito and T.~Sadahiro, \emph{Beta-expansions with negative bases}, Integers
  \textbf{9} (2009), A22, 239--259.

\bibitem{Kalle-Steiner}
C.~Kalle and W.~Steiner, \emph{Beta-expansions, natural extensions and multiple
  tilings associated with {P}isot units}, Trans. Amer. Math. Soc., to appear.

\bibitem{Klouda-Pelantova09}
K.~Klouda and E.~Pelantov{\'a}, \emph{Factor complexity of infinite words
  associated with non-simple {P}arry numbers}, Integers \textbf{9} (2009), A24,
  281--310.

\bibitem{Liao-Steiner}
L.~Liao and W.~Steiner, \emph{Dynamical properties of the negative
  beta-transformation}, arXiv:1101.2366v2.

\bibitem{Masakova-Pelantova}
Z.~Mas{\'a}kov{\'a} and E.~Pelantov{\'a}, \emph{{I}to-{S}adahiro numbers vs.
  {P}arry numbers}, arXiv:1010.6181v1.

\bibitem{Parry60}
W.~Parry, \emph{On the {$\beta$}-expansions of real numbers}, Acta Math. Acad.
  Sci. Hungar. \textbf{11} (1960), 401--416.

\bibitem{Renyi57}
A.~R{\'e}nyi, \emph{Representations for real numbers and their ergodic
  properties}, Acta Math. Acad. Sci. Hungar. \textbf{8} (1957), 477--493.

\bibitem{Thurston89}
W.~Thurston, \emph{Groups, tilings and finite state automata}, AMS Colloquium
  lectures (1989).

\end{thebibliography}
\end{document}